\documentclass[fontsize=12pt,a4paper,headings=normal,
twoside=false,leqno,parskip=half-,abstract=true]{scrartcl}
\usepackage[english]{babel}
\usepackage[utf8]{inputenc}
\usepackage{hyperref}
\hypersetup{
 pdftitle={},
 pdfauthor={},
 colorlinks=true,
 linkcolor=blue,
 citecolor=blue,
 filecolor=blue,
 urlcolor=blue}

\usepackage[pagewise]{lineno}

\usepackage[format=plain,labelfont=bf,font=small]{caption}
\usepackage{subfigure}
\usepackage{xcolor}
\usepackage[arrow, matrix, curve]{xy}
\usepackage{tikz}
\usepackage{float}

\usepackage{tikz-cd}
\usepackage{caption}
\usepackage{amsmath,amsthm}
\usepackage{amssymb} 
\usepackage[normalem]{ulem}

\newtheorem{theorem}{Theorem}[section]
\newtheorem{lemma}[theorem]{Lemma}
\newtheorem{cor}[theorem]{Corollary}

\theoremstyle{definition}
\newtheorem{definition}[theorem]{Definition}

\theoremstyle{remark}
\newtheorem{remark}[theorem]{Remark}

\numberwithin{equation}{section}

\DeclareMathAlphabet{\mathpzc}{OT1}{pzc}{m}{it}

\definecolor{refkey}{rgb}{1,0,0}
\definecolor{labelkey}{rgb}{1,0,0}

\title{Symbolic hunt of instabilities and bifurcations in reaction networks}
\author{Nicola Vassena }
\date{\today}

\begin{document}

\maketitle

\begin{abstract}
    The localization of bifurcations in large parametric systems is  still a challenge where the combination of rigorous criteria and informal intuition is often needed. With this motivation, we address symbolically the Jacobian matrix of reaction networks with general kinetics. More specifically, we consider any nonzero partial derivative of a reaction rate as a free positive symbol. The main tool are the Child-Selections: injective maps that associate {to a species $m$ a reaction $j$ where $m$ participates as reactant}. Firstly, we employ a Cauchy-Binet analysis and we structurally express any coefficient of the characteristic polynomial of the Jacobian in terms of Child-Selections. In particular, we fully characterize sign-changes of any of the coefficients. Secondly, we prove that the (in)stability of the Jacobian is inherited from the (in)stability of simpler submatrices identified by the Child-Selections. Thirdly, we provide sufficient conditions for purely imaginary eigenvalues of the Jacobian, hinting at Hopf bifurcation and oscillatory behavior.  All conditions are in terms of signs of integer stoichiometric submatrices identified by the Child-Selections {and do not require any Hurwitz-type computaton}. \mbox{Finally}, we focus on systems endowed with Michaelis-Menten kinetics and we show that any symbolic realization of the Jacobian matrix can be achieved at a fixed equilibrium by a proper choice of the kinetic constants.
\end{abstract}
\tableofcontents

\section{Introduction}

A major obstacle in the analysis of biological, ecological, and epidemiological networks is the uncertainty of the numerical parameters in the system. For this reason, dynamical models are typically given in general parametric forms, comprising many different and independent parameters. Even for the most advanced continuation software, it is at present still mandatory to fix most of the parameters with a certain degree of arbitrariness. Thus, it is still crucial for numerical analysis to develop some intuition about which parameter areas may show interesting dynamical behavior. The present paper is a contribution in this direction.

Let us focus on bifurcation behavior: we refer to the book by Guckenheimer and Holmes \cite{GuHo84} for general background. Vaguely speaking, a bifurcation is a sudden qualitative change in the system behavior according to a small change in the parameter values. One reason for interest in finding bifurcations is that they constitute the gates to parameter areas where interesting dynamical behavior occurs. For example, a \emph{saddle-node} bifurcation indicates a nearby area of multistationarity, while \emph{Hopf} bifurcations mark the birth of oscillations in form of periodic orbits. Multistationarity is the coexistence of multiple equilibria under otherwise identical conditions, and it is believed at the core of many epigenetic processes, including cell differentiation \cite{ThomKauf01}. Oscillations in biochemical systems are central in the regulation of metabolic processes, circadian rhythms, and other important biological functions \cite{Hess71}. The aforementioned bifurcations are local, that is, they are implied by certain conditions at one equilibrium. This suggests that finding a bifurcation point is often easier than directly proving the global dynamical behavior for which the bifurcation stands as a gate in the parameter space. To be perhaps more explicit, finding Hopf bifurcations is a standard method to infer oscillatory behavior in a system, albeit not all oscillations necessarily arise this way. 

We give a brief literature review of bifurcation analysis for biochemical networks.
To keep it concise, we focus on saddle-node and Hopf bifurcations. Both such bifurcations unfold by varying one single parameter, and they are hence the simplest and most studied cases. In the context of mass-action systems, saddle-node bifurcations have been addressed by Conradi et al. \cite{Conradi2007} and Domijan and Kirkilionis \cite{DomKirk09}. Here, the authors translated the general conditions for a saddle-node bifurcation in the polynomial language of mass action. Structural conditions for a saddle-node bifurcation in systems with general kinetics have been obtained in \cite{V23}, in the same setting as the present contribution. Gatermann et al. \cite{Gat2005} applied pioneering concepts from computer algebra to investigate Hopf bifurcations in mass-action systems. {Errami et al. developed a computational framework for Hopf bifurcation using convex coordinates \cite{Errami2015}}. Fiedler \cite{F19} proved global Hopf bifurcations in networks with feedback cycles that satisfy certain stoichiometric conditions. Through Hopf bifurcations in given biochemical systems, Conradi et al. detected oscillations in a mixed-mechanism phosphorylation system \cite{conradietal19}, Boros and Hofbauer in planar deficiency-one mass-action systems \cite{Boros21}, Hell and Rendall in the MAP kinase cascade \cite{Hell16}. Banaji and Boros \cite{BaBo23} fully classified all three-species, four-reactions bimolecular mass action systems that undergo Hopf bifurcation. Quite interesting and mathematically challenging is also the question of \emph{excluding} Hopf bifurcations from given systems. In this direction, structural necessary conditions have been addressed by Angeli et al. \cite{AnBaPa13}, while Conradi et al. \cite{CarstenHopfExclusion19} conjectured the absence of Hopf bifurcation in the sequential and distributive double phosphorylation cycle, studying few meaningful simplifications of the network. Hopf exclusion for the unsimplified network is still open, at present. {In the cited literature, the majority of the results employs the notorious Hurwitz criterion to find the spectral condition of a Hopf bifurcation for a choice of parameters, see \cite{Liu94} for an ad-hoc explanation. A major drawback of the Hurwitz approach is the great complexity of the computation: this approach becomes unfeasible for networks with more than very few species. In contrast, the results of the present paper do not argue via the Hurwitz criterion.}

Consider a dynamical system of ordinary differential equations:
\begin{equation}\label{firsteq}
    \dot{x}=g(x).
\end{equation}
An equilibrium bifurcation typically occurs when the Jacobian $G=\partial{g}/\partial{x}$ is nonhyperbolic, at an equilibrium $\bar{x}$. Here, nonhyperbolic means that at least one of its eigenvalues has zero-real part. Hence, the first step to finding and classifying bifurcations is studying the spectral properties of $G$. However, spectral conditions of $G$ are only necessary for a bifurcation to happen, but they are not sufficient to determine with precision the dynamical behavior nearby. Nevertheless, in applications, one might expect a certain bifurcation solely based on the spectral condition. For example, the class of vector fields undergoing a saddle-node bifurcation is open and dense in the class of smooth vector fields with an equilibrium with a singular Jacobian. Analogously, the class of vector fields undergoing a Hopf bifurcation is open and dense in the set of smooth vector fields with an equilibrium with a Jacobian with a pair of purely imaginary eigenvalues. Parameter areas are intrinsically uncertain and thus only qualitative behaviors that appear in open regions are meaningful. Thus, one might expect only such types of generic bifurcations. As a caveat, note that any application restricts the formulation of the problem and such a restriction may frustrate the generic expectation. For instance, the presence of a trivial equilibrium from which bifurcation occurs is a common assumption in many applications. For systems with such a structure, saddle-node bifurcations cannot happen, and \emph{transcritical} or \emph{pitchfork} bifurcations are likely to occur at an equilibrium with a singular Jacobian, depending on the symmetry of the system. In a biochemical setting, the network structure or the chosen class of reaction functions may as well obstruct the generic expectation. See for instance \cite{V22} for a network example where a singular Jacobian has always an algebraically double eigenvalue zero for any choice of reaction functions. In conclusion, the spectral approach can rule out a  bifurcation and it gives a solid guess for the occurrence of a certain bifurcation. Nevertheless, even the most solid guess must be checked afterward, either analytically or numerically.

This paper builds on this point of view and addresses symbolically the spectral properties of the Jacobian matrix $G$. We are interested in dynamical systems arising from reaction networks, or more general interaction networks of populations from ecology or epidemiology. Words, more than math, differ among such types of systems. For consistency, we mostly focus on and refer to chemical reaction networks. We consider the time-evolution $x(t)$ of the concentrations of the species or chemicals. Then, \eqref{firsteq} takes the form of 
$$\dot{x}=g(x):=Sf(x)$$
where $S$ is the stoichiometric matrix, in essence the incidence matrix of the network; $f(x)$ is the vector of the reaction functions. The precise form of $f$ is typically unknown in applications, and wide classes of parametric functions are used. For this reason, we consider $f$ from a quite general class of functions which we properly define in Definition \ref{monchemfunct} as \emph{monotone chemical functions}. In brief, we assume that a reaction function $f_j$ of a reaction $j$ is a monotone positive function only of the concentrations $x_m$ of the reactants $m$ of $j$. As addressed in more detail in Section \ref{sec_symbolic}, we interpret the Jacobian $G$ as a symbolic matrix, where the symbols are the partial derivatives of the reaction functions. At first, we study the spectral properties of $G$ independently from the dynamics. That is, we analyze how the network structures translate into the Jacobian, with no reference to any equilibrium. We discuss the \emph{realizability} of a spectral property at an equilibrium $\bar{x}$ for a proper choice of $f$ in a second moment. Such realizability is always possible within the entire class of monotone chemical functions: we simply think of a Taylor-type expansion to define functions with prescribed equilibrium values and derivatives. For applications that restrict the class of reaction functions to a given class of kinetics, it depends of course on the choice of kinetics. Mass action kinetics \cite{HFJ72}, in particular, does not provide enough parametric freedom for such an approach to be always valid, but exclusion results are still in full validity. For the likewise important class of enzymatic Michaelis-Menten kinetics \cite{MM13}, we already get that any spectral condition of the symbolic Jacobian can be realized at an equilibrium. We clarify this in Section \ref{MMsec}. Thus, the techniques presented in this paper provide equilibria with prescribed spectral properties for systems endowed with Michaelis-Menten kinetics.

More in detail, we have three main results.  The main tool for the analysis are Child-Selections, introduced in detail in Section \ref{sec_childselection}. A $n$-Child-Selection $\mathbf{J}^{(n)}$ is an injective map between $n$ species and $n$ reactions, such that $m$ is a reactant in the reaction $\mathbf{J}^{(n)}(m)$. Firstly, we employ the Cauchy-Binet formula to analyze the characteristic polynomial $\mathpzc{g}$ of the Jacobian $G$, Theorem \ref{CBTHM}. Such theorem extends to network settings the \emph{generalized Cauchy-Binet formula} \cite{Knill14}, proved by Knill. Moreover, we give structural network characterizations of the possibility of sign changes in any of the coefficients of $\mathpzc{g}$, Corollary \ref{corinst2}. Secondly, any $n$-Child-Selection $\mathbf{J}^{(n)}$ identifies a \emph{Cauchy-Binet component}, i.e. a $n\times n$ symbolic matrix of simple structure. In Theorem \ref{stability} we show that if any of such Cauchy-Binet components possesses $\sigma^+_{\mathbf{J}^{(n)}}>0$ unstable eigenvalues with positive-real part for some choice of symbols, then such instability is inherited by the full Jacobian for some choice of symbols. That is, the number $\sigma^+_G$ of unstable eigenvalues of the Jacobian $G$ is bigger or equal $\sigma^+_{\mathbf{J}^{(n)}}$. This allows the identification of easy and recognizable network motifs that lead to instability. As an example, Corollary \ref{autocat} states that an \emph{autocatalytic reaction}, where one of the reactants catalyzes a higher production of itself, is a source of instability. Lastly, in Lemma \ref{piprop} and Theorem \ref{pithm} we give sufficient conditions for the Jacobian to possess purely imaginary eigenvalues, hinting at Hopf bifurcations.

Besides the major inspiration from nonlinear dynamics, this paper can also be seen as a discussion on the stability properties of a certain class of symbolic matrices. Symbolic matrices have been addressed in the linear algebra literature \cite{BruSha09}, typically in form of \emph{sign patterns} $A$ with only prescribed sign of the entries from $\{+,-,0\}$. The Quirk-Ruppert-Maybee theorem characterizes the stability of such sign patterns, for any choice of the magnitude of the entries. See the focused account \cite{JKD77}. More recently, small-dimension sign patterns that exclude change of stabilities and bifurcations have been studied \cite{LOD18}. However, the symbolic class of sign patterns is too general for our purposes as it lacks the network structure this paper focuses upon. On another side, the structure of Jacobian matrices arising from networks often relates to studied classes of matrices that possess strong stability properties. For example, real $P_0$ matrices are defined by having all nonnegative principal minors. See the review by Hershkowitz \cite{Hersh92} in the linear algebra community, and Banaji et al. \cite{Ba-07} in a chemical network setting. In full analogy but more consistently for stability analysis of dynamical systems, $G$ is a $P^{(-)}_0$ matrix if $-G$ is a $P_0$ matrix. In particular, real positive eigenvalues are excluded for $P^{(-)}_0$ matrices, as a straightforward consequence of Descartes' rule of sign applied to their characteristic polynomial. However, $P^{(-)}_0$ matrices can still be unstable due to pairs of complex-conjugate eigenvalues with positive real part. Consider now any symbolic family $G(\boldsymbol{\rho})$ of $P^{(-)}_0$ matrices, where $\boldsymbol{\rho}$ indicates the dependence on one or more symbols. Assume there is a stable matrix $G(\bar{\boldsymbol{\rho}})$ in the family. Then, either the family is stable for all choices of $\boldsymbol{\rho}$, or there is a choice $\boldsymbol{\rho}^*$ such that $G(\boldsymbol{\rho}^*)$ has purely imaginary eigenvalues. In the case where $G(\boldsymbol{\rho})$ is a Jacobian of a dynamical system, saddle-node bifurcations are always excluded and only Hopf bifurcations are possible. In our setting, we address connections to $P^{(-)}_0$ matrices in Corollaries \ref{pmatrix}, \ref{exclusive}, and \ref{corlem}.

The paper is organized as follows: Sections \ref{sec_reactnet} and \ref{sec_symbolic} introduce reaction networks and the symbolic approach, respectively. We present in Section \ref{sec_childselection} the central tool: Child-Selections. The main results are collected in Section \ref{main}. More specifically, subsection \ref{CB} focuses on the Cauchy-Binet analysis for the characteristic polynomial of the Jacobian, subsection \ref{instabilities} discusses the inheritance of stability properties of Child-Selections to the full network, and subsection \ref{purelyim} provides sufficient conditions for the Jacobian to have purely imaginary eigenvalues. Section \ref{MMsec} establishes the validity of our symbolic approach for systems endowed with Michaelis-Menten kinetics. Section \ref{examplesec} provides two toy models where we exemplify our results. Example I explicitly finds Michaelis-Menten constants such that the system has an equilibrium with purely imaginary eigenvalues of the Jacobian. Example II uses the results to find two distinct parameter areas where multistationarity and oscillations occur: numerical simulations are provided. Section \ref{discussion} concludes the paper with the discussion section. The proofs of the main results are postponed in Section \ref{proofs}.

\textbf{Acknowledgments.} I thank Carsten Conradi for many inspiring discussions. This work has been supported by the DFG (German Research Foundation), project no. 512355535.

\section{Reaction networks}\label{sec_reactnet}

A \emph{chemical reaction network} is a set $\mathbf{M}$ of species or chemicals, together with a set $\mathbf{E}$ of reactions. The cardinalities of such sets is $|\mathbf{M}|=M$ and $|\mathbf{E}|=E$. Letters $m \in \mathbf{M}$ and $j \in \mathbf{E}$ refer to species and reactions, respectively. By labeling the network, we arbitrarily fix an order to the species set $\mathbf{M}=\{m_1,..,m_M\}$ and to the reaction set $\mathbf{E}=\{j_1,...,j_E\}$. In the examples, we use capital letters $A,B,C,...$ for species and natural numbers $1,2,3,..$ for reactions.

A reaction $j \in \mathbf{E}$ is an ordered association between nonnegative linear combinations of the species:
\begin{equation}\label{reactionj}
s^j_{m_1}m_1+...+s^j_{m_M}m_M \quad \underset{j}{\rightarrow}\quad \tilde{s}^j_{m_1}m_1+...+\tilde{s}^j_{m_M}m_M,
\end{equation}
with integer \emph{stoichiometric} coefficients $s^j_m,\tilde{s}^j_m\in \mathbb{Z}_{\ge 0}$. We could analogously consider real stoichiometric coefficients with no mathematical difference. The species $m$ appearing at the left side of \eqref{reactionj} with nonzero stoichiometric coefficient $s^j_m$ are called \emph{reactants} of the reaction $j$. Respectively, the species $m$ appearing at the right side of \eqref{reactionj} with nonzero stoichiometric coefficient $\tilde{s}^j_m$ are called \emph{products} of the reaction $j$. Frequently reaction networks are \emph{open system}, i.e., they exchange chemicals with the outside environment. For this reason, we also consider reactions with no outputs (\emph{outflow reactions}) or with no inputs (\emph{inflow reactions}).

We call a reaction $j_{aut}$ \emph{autocatalytic} in the chemical $m$ if 
 $$\tilde{s}^{j_{aut}}_m > s^{j_{aut}}_m > 0.$$
In particular, $m$ is both a reactant and a product of $j_{aut}$. An \emph{autocatalytic network} is simply a network that possesses at least one autocatalytic reaction.

Throughout, the notation $A^k_h$ indicates the entry in the $h^{th}$ row and $k^{th}$ column of a matrix $A$. The $M \times E$ \emph{stoichiometric matrix} $S$ is then defined as
\begin{equation}
S^j_m:=\tilde{s}^j_m-s^j_m.
\end{equation}
 This way we fix a direction to the reactions, considering them as \emph{irreversible}.  Reversible reactions $j$ such as
$$A\quad\underset{j}{\rightleftharpoons} \quad B$$
can be taken into account in our setting by considering two different reactions 
\begin{equation}\label{j1}
A\quad\underset{j_1}{\rightarrow} \quad B \quad \text{and}\quad B\quad\underset{j_2}{\rightarrow} \quad A.
\end{equation}

The notation $S^j$ indicates the $j^{th}$ column of $S$, corresponding to reaction $j$. For example, in a network of three species $A,B,C$, the stoichiometric column of reaction $j_1$ in \eqref{j1} reads:
$$S^{j_1}=\begin{pmatrix}
    -1\\
    1\\
    0
\end{pmatrix}.$$

We use the variable $x$ to indicate the $M$-vector of the concentrations of the species. We consider strictly positive concentrations $x>0$ and we proceed under the assumption that the reactor is well-mixed, spatially homogeneous, and isothermal. The time-evolution $x(t)$ satisfies the system of ODEs:
\begin{equation}\label{mainsyst}
    \dot{x}=Sf(x),
\end{equation}
where $S$ is the $M\times E$ stoichiometric matrix and $f(x)$ is the $E$-vector of the reaction functions. With no reactant, we consider the reaction function of inflow reactions $j_F$ as constant:
$$f_{j_F} \equiv F.$$
For any other reaction $j$, we assume that $f_j$ is monotone chemical, defined as follows.
\begin{definition}[Monotone chemical functions]\label{monchemfunct}
Let $j$ be a reaction and $f_j$ the associated reaction function. We call $f_j$ \emph{chemical} if 
\begin{enumerate}
\item $f_j$ depends only on the concentrations of the reactants of the reaction $j$;
\item $f_j$ is positive, i.e.,
$$f(x)>0,\text{ for every $x>0$}.$$
\end{enumerate}
We call a chemical function $f_j$ \emph{monotone} if
\begin{enumerate}
\setcounter{enumi}{2}
\item $f'_{jm}(x):=\frac{\partial f_j (x)}{\partial x_m} > 0$,  for any reactant species $m$ of $j$ and $x>0$.
\end{enumerate}
\end{definition}
Widely used kinetics as mass action, Michaelis--Menten, and Hill kinetics \cite{Hill10} follow Definition \ref{monchemfunct}. We highlight however few restrictions. Condition 1 excludes dependencies $f'_{jm}(x) \neq 0$ not expressed by the stoichiometry. Regulatory terms, i.e. $f'_{jm}(x) \neq 0$ with $m$ not a reactant to $j$, both in form of \emph{activators} $f'_{jm}(x) > 0$ and \emph{inhibitors} $f'_{jm}(x) < 0$, are not taken in account. Further, we do not consider nonmonotone reaction rates such as, e.g., substrate inhibition. On the other hand, monotone \emph{decreasing} functions, i.e. $f'_{jm}(x) < 0$, could be considered analogously, carrying no essential difference. We have chosen monotone increasing functions as this case appeared more relevant for applications: the decreasing case can be easily adapted from the presented results.

\section{A symbolic approach} \label{sec_symbolic}

We aim at studying stability properties of equilibria $\bar{x}$ of \eqref{mainsyst}:
\begin{equation}\label{equilibria}
    0=Sf(\bar{x}).
\end{equation}
We only consider networks whose stoichiometric matrix $S$ admits a positive right kernel vector $\mathbf{r}$
\begin{equation}\label{equilibriaflux}
S\mathbf{r}=0,
\end{equation}
with $r_j > 0$ for all reactions $j$. Without this assumption, there would be no choice of monotone chemical functions $f$ that admits an equilibrium. On the other hand, this assumption guarantees the existence of choices of reaction functions $f$ satisfying \eqref{equilibria}, at a given $\bar{x}$. {Networks, whose stoichiometric matrix satisfies \eqref{equilibriaflux} for a positive vector $\mathbf{r}$, are called \emph{consistent} in the literature \cite{Ang07} or also \emph{dynamically nontrivial} \cite{BaBo23}}. The stability of an equilibrium $\bar{x}$ can be discussed at first approximation by studying the Jacobian $G|_{\bar{x}}$ of \eqref{mainsyst}. In particular, the matrix stability of a \emph{hyperbolic} $G|_{\bar{x}}$, i.e. without zero-real-part eigenvalue, is inherited by the 
dynamical stability of the equilibrium $\bar{x}$. Complementarily, the case when $G|_{\bar{x}}$ is not hyperbolic is particularly interesting as it points to bifurcations and qualitative changes in the dynamics.

We study the Jacobian $G$ symbolically. That is, we consider $G$ as a symbolic matrix where the positive symbols $r'_{jm}$ correspond to the partial derivatives $f'_{jm}$. We use different letters, $r'_{jm}$ vs $f'_{jm}$, precisely to stress that the symbols $r'_{jm}$ are free positive parameters and they are independent - at first - from a specific choice of reaction functions $f$. The bold $\mathbf{r}'$ indicates the set of such positive symbols and $G(\bar{\mathbf{r}}')$ refers to the symbolic Jacobian $G$ evaluated at the choice $\bar{\mathbf{r}}'$. Throughout the paper, we often deal with submatrices and minors of $G$. To avoid an overload of notation, we do not specify further whether or not all the symbols are present in the considered matrix. For example, two choices of symbols for the first diagonal entry of $G$ are again denoted by 
$G_1^1(\bar{\mathbf{r}}'_1), G_1^1(\bar{\mathbf{r}}'_2)$ even if possibly not all symbols $\mathbf{r}'$ appear in  $G^1_1$.

The $E \times M$ \emph{reactivity matrix} of the symbolic partial derivatives is defined as follows.
\begin{equation}
    R_j^m:=
    \begin{cases}
        r'_{jm} \text{ if $f'_{jm}\neq 0$};\\
        0 \quad \text{otherwise}.
    \end{cases}
\end{equation}
This way the symbolic Jacobian $G$ can be expressed as a product of two matrices: the integer-valued stoichiometric matrix $S$ and the symbolic-valued reactivity matrix $R$:
$$G:=SR.$$

The question of whether a certain choice of symbols $\mathbf{r'}$ can be realized at a given equilibrium $\bar{x}$ and for a given choice of parametric functions $f$ depends of course on the parametric freedom of the chosen class of reaction functions $f$ (kinetics). In the wide class of chemical functions defined in Definition \ref{monchemfunct}, it is clear that we can always find reaction functions $f$ satisfying at a chosen value $\bar{x}$
\begin{equation}
    \begin{cases}
        f_j(\bar{x})= r_j\\
        f'_{jm}(\bar{x})=r'_{jm}
    \end{cases},
\end{equation} 
for any reaction $j$ and species $m$. We just think of a Taylor-type expansion to define the desired function. However, applications typically restrict the parametric functions $f$ endowing the network. For example, \emph{mass action kinetics} reads
\begin{equation}\label{ma}
    f_j (x) := k_j \prod_{m\in \mathbf{M}} x_m^{s^j_m},
\end{equation}
where $k_j >0$ is a positive constant and the integer $s^j_m$ is the stoichiometric coefficient of the species $m$ as reactant of the reaction $j$. The value of any derivative
\begin{equation} \label{mader}
f'_{jm}(x)=\;{s^j_m}\; x_m^{(s^j_m - 1)} \; k_j \prod_{n \neq m}  x_n^{s^j_n}=\frac{s^j_m}{x_m} f_j
\end{equation}
cannot be chosen independently from the function value $f_j(\bar{x})$ at a fixed $\bar{x}$. This partially excludes mass action kinetics for the current symbolic approach. More precisely,
\emph{exclusion results} that forbid a spectral property of $G$ for any choice of monotone chemical functions are valid also in the mass action case. On the contrary, however, \emph{existence results} of certain spectral properties of $G$ do not conclude the realizability of the same spectral properties in the system endowed with mass action kinetics. Slightly richer kinetics already guarantee the realizability of any $\mathbf{r}'$ at a fixed equilibrium $\bar{x}$. We address explicitly the case of Michaelis-Menten kinetics in Section \ref{MMsec}. This serves both as a specific example and as a guidance on how to implement the results in given dynamical models.

\section{$n$-Child-Selections}\label{sec_childselection}

For $n\le M$, a set $\mathbf{M}^{(n)}$ consists
of any choice of $n$ chemicals out of the set $\mathbf{M}$. We consider the order of any $\mathbf{M}^{(n)}$ as induced by the order of $\mathbf{M}$. To avoid overload of notation, we do not assign an index to this family of sets, but we only use $\bar{n}$ or $\bar{\mathbf{M}}^{(n)}$ to refer to a fixed choice of $n$ and  $\mathbf{M}^{(n)}$, respectively. 
 
\begin{definition}[$n$-Child-Selection]
Consider any set $\bar{\mathbf{M}}^{(n)}$. A $n$-Child-Selection $\mathbf{J}^{(n)}$ is an injective map
$$\mathbf{J}^{(n)}: \quad \bar{\mathbf{M}}^{(n)} \quad \mapsto \quad \mathbf{E},$$
associating to each chemical $m \in \bar{\mathbf{M}}^{(n)}$ a reaction $j$ such that $m$ is a reactant to $j$.
\end{definition}
In particular, note that
$$r'_{\mathbf{J}^{(n)}(m)m}\neq 0,\text{ for $m$ in $\bar{\mathbf{M}}^{(n)}$.}$$
To simplify as possible the notation, we refer to Child-Selections $\mathbf{J}^{(n)}$ without an explicit reference to its domain $\bar{\mathbf{M}}^{(n)}$. 

An $n$-Child-Selection uniquely identifies an $n \times n$ integer matrix $S[\mathbf{J}^{(n)}]$ constructed from the stoichiometric matrix as follows. Let $m_1, ..., m_n$ be the $n$ chemicals in the set $\bar{\mathbf{M}}^{(n)}$. Then $S^{\mathbf{J}^{(n)}}$ indicates the $M\times n$ matrix such that the $i^{th}$ column is the stoichiometric column $S^j$ of the reaction $j=\mathbf{J}^{(n)}(m_i)$, and \begin{equation}\label{Sjsubmatrix}
S[\mathbf{J}^{(n)}]:=S^{\mathbf{J}^{(n)}}_{\bar{\mathbf{M}}^{(n)}} 
\end{equation}
is the $n\times n$ square matrix obtained from $S^{\mathbf{J}^{(n)}}$ by removing the rows corresponding to chemicals $m \not\in {\bar{\mathbf{M}}^{(n)}}$. In essence, $S[\mathbf{J}^{(n)}]$ is a reshuffled square submatrix of the stoichiometric matrix $S$.

Moreover, we associate to any Child-Selection a \emph{behavior coefficient}: 
\begin{equation}\label{behaviorcoefficient}
    \alpha_{\mathbf{J}^{(n)}}:= \operatorname{det} S[\mathbf{J}^{(n)}].
\end{equation}

Based on previous work of the author \cite{VGB20}, it is possible to structurally characterize the behavior coefficient by analyzing certain cycles in the network. 

Finally, $R[\mathbf{J}^{(n)}]$ indicates  the $n\times n$ symbolic diagonal matrix defined as:
$$R[\mathbf{J}^{(n)}]_m^m:=r'_{\mathbf{J}^{(n)}(m)m}, \quad \text{for $m \in \bar{\mathbf{M}}^{(n)}$.}$$

\section{Main results}\label{main}

\subsection{Cauchy-Binet expansion for the characteristic polynomial}\label{CB}

We start with a definition.

\begin{definition}[Cauchy-Binet component]\label{defcb}
Consider an $n$-Child-Selection $\mathbf{J}^{(n)}$. We call the symbolic $n \times n$ matrix
$$G[\mathbf{J}^{(n)}]:=S[\mathbf{J}^{(n)}]R[\mathbf{J}^{(n)}]$$
the \emph{Cauchy-Binet component} of $G$ associated to $\mathbf{J}^{(n)}$.
\end{definition}

\begin{remark}\label{dstability}
A Cauchy-Binet component is thus defined as the product of a matrix $S[\mathbf{J}^{(n)}]$ with a positive diagonal matrix $R[\mathbf{J}^{(n)}]$. Such algebraic structure has already been discussed in the literature. In fact, \emph{$D$-stable (resp. $D$-hyperbolic) matrices} are defined as matrices $A$ such that $B:=AD$ is stable (resp. hyperbolic) for any choice of positive diagonal matrix $D$. See \cite{Gio15,Ku19} for an overview. The problem of characterizing the class of $D$-stable matrices $A$ is open.
\end{remark}

\begin{remark}\label{remarksummands}
A Cauchy-Binet component $G[\mathbf{J}^{(n)}]$ can also be interpreted as the symbolic matrix obtained by multiplying each column of $S[\mathbf{J}^{(n)}]$ for a symbol. For $\mathbf{J}^{(n)}: \mathbf{M}^{(n)}\mapsto \mathbf{E}$, in particular, we have 
\begin{equation}\label{cbmonomial}
 \operatorname{det}G[\mathbf{J}^{(n)}]=\operatorname{det}S[\mathbf{J}^{(n)}]\operatorname{det}R[\mathbf{J}^{(n)}]=\alpha_{\mathbf{J}^{(n)}} \prod_{m\in\mathbf{M}^{(n)}} r'_{\mathbf{J}^{(n)}(m)m},
\end{equation}
where $\alpha_{\mathbf{J}^{(n)}}$ is the behavior coefficient \eqref{behaviorcoefficient} of $\mathbf{J}^{(n)}$. In particular, $$\operatorname{sign}\operatorname{det}G[\mathbf{J}^{(n)}]=\operatorname{sign}\operatorname{det}S[\mathbf{J}^{(n)}].$$
\end{remark}
A central lemma follows:
\begin{lemma}\label{corinst}
Consider any set of $n$ species $\mathbf{M}^{(n)}$. Then, the principal minor $\operatorname{det}G^{\mathbf{M}^{(n)}}_{{\mathbf{M}}^{(n)}}$ of the Jacobian matrix $G$ can be expanded as
\begin{equation}\label{corinstsum}
\operatorname{det}G^{\mathbf{M}^{(n)}}_{{\mathbf{M}}^{(n)}}=\sum_{\mathbf{J}^{(n)}:\mathbf{M}^{(n)}\;\mapsto \mathbf{E}} \operatorname{det} G[\mathbf{J}^{(n)}].
\end{equation}
\end{lemma}
{The sum \eqref{corinstsum} runs over all $n$-Child-Selections defined on the set $\mathbf{M}^{(n)}$}. The proof of Lemma \ref{corinst} relies on the Cauchy-Binet formula. A direct consequence is the following theorem. 
\begin{theorem}\label{CBTHM}
Consider the Jacobian matrix $G$ of the system \eqref{mainsyst} and its characteristic polynomial
\begin{equation}\label{charpolform}
\mathpzc{g}(\lambda):=\operatorname{det}(G-\lambda \operatorname{Id})=(-\lambda)^M+a_1(-\lambda)^{M-1}+...+a_{n}(-\lambda)^{M-n}+...+a_M.
\end{equation}
Then, for each coefficient $a_n$, the following expansion holds
\begin{equation}\label{cbcharpol}
    a_n=\sum_{\mathbf{J}^{(n)}} \operatorname{det} G[\mathbf{J}^{(n)}].
\end{equation}
\end{theorem}
The sum \eqref{cbcharpol} runs over all $n$-Child-Selections. {We underline that each $n$-Child-Selection identifies a single monomial \eqref{cbmonomial}, and thus cannot be further expanded. In this sense, if we consider $a_n$ as a multilinear $n$-homogenous polynomials int the variables $\mathbf{r}'$, then \eqref{cbcharpol} uniquely expands the coefficient $a_n$ in a sum of monomials.} We proceed with a definition.
\begin{definition}[Fixed sign]\label{fixedsign}
    A symbolic expression $P(\mathbf{r}')$ in the positive symbols $\mathbf{r}'>0$ is said to be of \emph{fixed sign} if 
    its sign does not depend on the choice of the symbols $\mathbf{r}'$. 
\end{definition}

 We have the following corollary.

\begin{cor}\label{corinst2}
A principal minor $\operatorname{det}G^{\bar{\mathbf{M}}^{(n)}}_{\bar{\mathbf{M}}^{(n)}}$ is of fixed sign if and only if 
$$\alpha_{\mathbf{J}_1^{(n)}}\alpha_{{\mathbf{J}_2^{(n)}}} \ge 0,$$
for any two $n$-Child-Selections $\mathbf{J}_1^{(n)},\mathbf{J}_2^{(n)}$ such that
$$
\mathbf{J}_i^{(n)}:\bar{\mathbf{M}}^{(n)}\mapsto \mathbf{E}, \quad \text{for $i=1,2$.}$$
In general, the coefficient $a_n$ is of fixed sign if and only if
$$\alpha_{\mathbf{J}_1^{(n)}}\alpha_{{\mathbf{J}_2^{(n)}}} \ge 0,$$
for any two $n$-Child-Selections $\mathbf{J}_1^{(n)},\mathbf{J}_2^{(n)}$.
\end{cor}

We conclude this subsection with a corollary about $P^{(-)}_0$ matrices. We recall that a matrix $A$ is $P^{(-)}_0$ matrix if all of its nonzero principal minors of order $n$ have sign $(-1)^n$. We repeat that nonzero Jacobians that are also $P^{(-)}_0$ matrices possess interesting properties for the dynamics. The most important is that multistationarity is excluded. In particular, saddle-node bifurcations are excluded and no real positive eigenvalues are possible. However, an equilibrium whose Jacobian is a $P^{(-)}_0$ matrix may nevertheless be unstable due to the presence of pairs of complex-conjugate eigenvalues with positive real part. In \cite{Ba-07}, the authors presented a characterization of reaction networks that possess $P^{(-)}_0$ Jacobians for any choice of parameters. In our slightly different setting and language, Corollary \ref{corinst2} implies a straightforward characterization of Jacobians that are  $P^{(-)}_0$ for all choices of symbols $\mathbf{r}'$.
\begin{cor}[$P^{(-)}_0$ matrix I]\label{pmatrix}
The Jacobian $G$ is a $P^{(-)}_0$ matrix for all choices of symbols $\mathbf{r}'$ if and only if for all Child-Selections $\mathbf{J}^{(n)}$ it holds
$$\operatorname{sign}\alpha_{\mathbf{J}^{(n)}}=(-1)^{n}.$$
\end{cor}

\subsection{Instabilities}\label{instabilities}

We call \emph{stable} (\emph{unstable}) an eigenvalue with negative (positive) real part. An $M\times M$ matrix with $M$ stable eigenvalues is called a \emph{stable matrix}, whereas it is called \emph{unstable} if it has at least one unstable eigenvalue. Accordingly, the expression \emph{``the network admits stability''} refers to the existence of a choice of symbols $\mathbf{r}'$ such that the Jacobian $G(\mathbf{r}')$ evaluated at $\mathbf{r}'$ is a stable matrix. Complementarily, the expression \emph{``the network admits instability''} means that there exists a choice of symbols $\mathbf{r}'$ such that the Jacobian $G(\mathbf{r}')$ evaluated at $\mathbf{r}'$ is unstable. We recall the concept of \emph{inertia} of a square matrix. 
\begin{definition}[Inertia of a matrix]
    The inertia of an $M \times M$ square matrix $A$ is a nonnegative triple  $$(\sigma^-_A,\sigma^+_A,\sigma^0_A),$$
where $\sigma^-_A$,$\sigma^+_A$ and $\sigma^0_A$
 are the number of stable, unstable, and zero-real-part eigenvalues of $A$, respectively. The eigenvalues are counted with their multiplicities so that $\sigma^+_A+\sigma^-_A+\sigma^0_A=M$.
 \end{definition}

 The main result of this section follows.
 
\begin{theorem}\label{stability}
    Consider any $n$-Child-Selection $\mathbf{J}^{(n)}:\mathbf{M}^{(n)} \mapsto \mathbf{E}$ and its associated Cauchy-Binet component $G[\mathbf{J}^{(n)}]$. Assume there exists a choice $\mathbf{r}'_1$ of symbols such that
    $$\operatorname{inertia}(G[\mathbf{J}^{(n)}](\mathbf{r}'_1))=(\sigma^-_{\mathbf{J}^{(n)}},\sigma^+_{\mathbf{J}^{(n)}},\sigma^0_{\mathbf{J}^{(n)}}).$$
    Then there exists a choice of $\mathbf{r}_2'$ such that
    $$\operatorname{inertia}(G(\mathbf{r}_2'))=(\sigma^-_G,\sigma^+_G,\sigma^0_G),$$
    with $\sigma^-_G \ge \sigma^-_{\mathbf{J}^{(n)}}$, and $\sigma^+_G \ge \sigma^+_{\mathbf{J}^{(n)}}$. 
\end{theorem}

 We derive four corollaries of interest. The first provides a sufficient condition for stability in the network.
\begin{cor}\label{stable}
    Recall the cardinality of the set of species $|\mathbf{M}|=M$. Assume that one $M$-Child-Selection $\mathbf{J}^{(M)}$ is such that the associated Cauchy-Binet component $G[\mathbf{J}^{(M)}]$ is stable for a choice of symbols. Then the network admits stability.\\ 
    In particular, if the associated integer matrix $S[\mathbf{J}^{(M)}]$ is stable, then the network admits stability.
\end{cor}
 The second provides a sufficient condition for instability in the network.
\begin{cor}\label{cor1}
Assume that one $n$-Child-Selection $\mathbf{J}^{(n)}$ has coefficient behavior s.t. 
\begin{equation}\label{badsign}
    \operatorname{sign}\alpha_{\mathbf{J}^{(n)}}=(-1)^{n-1}.
\end{equation}
Then the network admits instability.
\end{cor}

Corollary \ref{cor1} constitutes the complement to Corollary \ref{pmatrix}: together they form the following.

\begin{cor}[$P^{(-)}_0$ matrix II]\label{exclusive}
If the Jacobian is not a $P^{(-)}_0$ matrix for some choice of symbols $\mathbf{r}'$, then the network admits instability. 
\end{cor}

In particular, if the network does not admit instability then the Jacobian is a $P^{(-)}_0$ matrix for all choices of symbols. It is anyways possible to have unstable $P^{(-)}_0$ matrices: a loss of stability can occur via purely imaginary crossings. The last corollary focuses on a common source of instability: autocatalysis. 

\begin{cor}[Autocatalysis I]\label{autocat}
Autocatalytic networks admit instability.
\end{cor}

See \cite{VasStad23} for an extensive analysis of the relation between autocatalysis and instability, in the setting of the present work.

\subsection{Purely imaginary eigenvalues}\label{purelyim}

We start with a lemma.
\begin{lemma}\label{piprop}
Let $\mathbf{J}^{(n)}$ be any Child-Selection. Assume there are two choices $\mathbf{r}_1'$, $\mathbf{r}_2'$ of symbols such that the associated Cauchy-Binet component $G[\mathbf{J}^{(n)}](\mathbf{r}')$ is stable for $\mathbf{r}'=\mathbf{r}'_1$ and unstable for $\mathbf{r}'=\mathbf{r}'_2$. Then there exists a choice $\mathbf{r}_3'$ of symbols such that the Jacobian $G(\mathbf{r}_3')$ has purely imaginary eigenvalues. 
\end{lemma}
\begin{remark}
We stress how Lemma \ref{piprop} as well as the following results of this section do not automatically imply further nonresonance conditions such as the absence of other zero-real-part eigenvalues. 
\end{remark}
We state a corollary of interest. 
\begin{cor}[$P^{(-)}_0$ matrix III]\label{corlem}
Let $\mathbf{J}^{(n)}:\bar{\mathbf{M}}^{(n)}\mapsto \mathbf{E}$ be a $n$-Child-Selection. Assume that the associated matrix $S[\mathbf{J}^{(n)}]$ is a stable matrix that is not a $P^{(-)}_0$ matrix. Then there exists a choice $\mathbf{r}'$ of symbols such that $G(\mathbf{r}')$ has purely imaginary eigenvalues. 
\end{cor}
\begin{remark}[Autocatalysis II] The assumptions of Corollary \ref{corlem} are met if $S[\mathbf{J}^{(n)}]$ is stable and there exists a species $\bar{m}\in \bar{\mathbf{M}}^{(n)}$ such that $\bar{j}=\mathbf{J}^{(n)}(\bar{m})$ is autocatalytic. 
\end{remark}

We proceed to the main result of this section.

\begin{theorem}\label{pithm}
Fix $\bar{n}\in \{1,...,M\}$ and a set $\bar{\mathbf{M}}^{(\bar{n})}$. 
Consider a collection $\mathpzc{C}$ of $\bar{n}$-Child-Selections $\{\mathbf{J}_1^{(\bar{n})},...,\mathbf{J}_k^{(\bar{n})}\}$ such that
\begin{equation}\label{invertibility}
\text{either $\operatorname{sign}\alpha_{\mathbf{J}^{(\bar{n})}}= (-1)^{\bar{n}}$ or $\operatorname{sign}\alpha_{\mathbf{J}^{(\bar{n})}}=0$}
\end{equation}
for all $\mathbf{J}^{(\bar{n})}\in \mathpzc{C}$. Call $\mathpzc{P}$ the set of pairs $(m,j)$ such that $m\in \bar{\mathbf{M}}^{(\bar{n})}$ and $j=\mathbf{J}^{(\bar{n})}(m)$ for some $\mathbf{J}^{(\bar{n})}\in \mathpzc{C}$. Assume the following two conditions:
\begin{enumerate}
    \item there exists a $\bar{n}$-Child-Selection $\mathbf{J}_1^{(\bar{n})} \in \mathpzc{C}$ such that $$S[\mathbf{J}_1^{(\bar{n})}]\text{ is a stable matrix.}$$
    \item there exists a $n$-Child-Selection $\mathbf{J}_2^{(n)}$, $n < \bar{n}$,
$$\mathbf{J}_2^{(n)}:\mathbf{M}^{(n)}\mapsto \mathbf{E},$$ with
$\mathbf{M}^{(n)}\subset \bar{\mathbf{M}}^{(\bar{n})}$ 
   and $(m,\mathbf{J}_2^{(n)}(m))\in\mathpzc{P}$ for all $m \in \mathbf{M}^{(n)}$, such that $$S[\mathbf{J}_2^{(n)}]\text{ is an unstable matrix.}$$
\end{enumerate}
Then, there exists a choice $\mathbf{r}'$ of symbols such that the Jacobian $G(\mathbf{r}')$ admits purely imaginary eigenvalues.
\end{theorem}

\begin{remark}
    The idea of the proof is as follows. Condition 1 implies the existence of a choice $\mathbf{r}_1$ of symbols such that $G[\mathbf{J}^{(\bar{n})}_1](\mathbf{r}'_1)$ is stable, while condition 2 implies the existence of a choice $\mathbf{r}_2$ of symbols such that $G[\mathbf{J}_2^{(n)}](\mathbf{r}'_2)$ is unstable. Under the invertibility assumption \eqref{invertibility} and relying on the intermediate value theorem we find purely imaginary eigenvalues. We can generalize Theorem \ref{pithm} by substituting conditions 1 and 2 with the bare existence of choices $\mathbf{r}'_1$ and $\mathbf{r}'_2$ of symbols such that the Cauchy-Binet components $G[\mathbf{J}^{(\bar{n})}_1](\mathbf{r}'_1)$ and $G[\mathbf{J}_2^{(n)}](\mathbf{r}'_2)$ are stable and unstable, respectively. The proof is analogous: Theorem \ref{pithm} is technically just a corollary of such a result. However, there are two reasons why Theorem \ref{pithm} is stated this way. Firstly, checking stability conditions for symbolic matrices is a formidable task. On the contrary, the conditions of Theorem \ref{pithm} rely only on integer matrices and they can be easily found even in reasonably large networks. In particular, condition 2 is already implied by $\operatorname{det}S[\mathbf{J}_2^{(n)}]=(-1)^{(n-1)}$. Secondly, assume there exists $\mathbf{r}_3'$ such that $G[\mathbf{J}_1^{(\bar{n})}](\mathbf{r}_3')$ is unstable, or respectively $G[\mathbf{J}_2^{(n)}](\mathbf{r}_3')$ is stable. Lemma \ref{piprop} already concludes there is a choice of symbols such that $G$ has purely imaginary eigenvalues.
\end{remark}

For the remarkable case in which the Jacobian $G$ is of fixed sign, we have the following straightforward corollary.
\begin{cor}\label{corpi}
Assume that the Jacobian $G$ has a nonzero determinant of fixed sign. If the network admits both stability and instability, then there exists a choice of $\mathbf{r}'$ such that $G(\mathbf{r}')$ has purely imaginary eigenvalues.
\end{cor}

\section{Realizability for Michaelis-Menten kinetics}\label{MMsec}

We show the applicability of the results to the reaction scheme known as Michaelis-Menten kinetics. We follow and expand an argument from \cite{V23}. Michaelis-Menten kinetics is typically used to model enzymatic reactions in cellular metabolism. Independently, the same mathematical form corresponds to \emph{Holling type II functional response} in ecological models \cite{Holling65}, and to the \emph{Monod equation} for the growth of microorganisms \cite{Monod49}. We proceed referring only to Michaelis-Menten. The mathematical form of a reaction $j$ according to Michaelis-Menten kinetics is
\begin{equation}\label{MMeq}
f_j(x):=a_j\prod_{m\in\mathbf{M}} \Bigg( \frac{x_m}{(1+b^j_m x_m)}\Bigg)^{s^j_m},
\end{equation}
where $a_j > 0$, $b^j_m \ge 0$ and the integer $s^j_m$ is the stoichiometric coefficient of species $m$ as reactant of the reaction $j$. Mass action kinetics is recovered by considering the case $b^j_m=0$ for all $j,m$. We use the bold $\mathbf{a}$ for the set of parameters $a_j$ for all reactions $j$. Respectively, we denote $\mathbf{b}$ the set of parameters $b^j_m$, for all $j$ and $m$. 
\begin{theorem}\label{MM}
 Fix any positive triple $(\bar{x},\bar{\mathbf{r}},\bar{\mathbf{r}}')>0$. That is, fix a positive concentration value $\bar{x}>0$, a positive vector of fluxes $\bar{\mathbf{r}}>0$, and a positive vector of partial derivatives symbols $\bar{\mathbf{r}}'>0$. Then there exists a choice of parameters $\mathbf{a},\mathbf{b}$ such that the Michaelis-Menten function $f_j(\bar{x})$ satisfies
 \begin{equation}\label{mmconstr}
     \begin{cases}
     f_j(\bar{x})= K \bar{r}_j;\\
     \partial f_j (\bar{x})/ \partial x_m = \bar{r}'_{jm},
     \end{cases}
 \end{equation}
for any reaction $j$ with reactants $m$. The positive constant $K>0$ is independent of $j$. In particular, if  $\bar{\mathbf{r}}$ is an equilibrium flux, i.e. $S\bar{\mathbf{r}}=0$, then the concentration $\bar{x}$ is a positive equilibrium of the Michaelis-Menten system
$$\dot{x}=S f(x),$$
with prescribed partial derivatives $\bar{\mathbf{r}}'$.
\end{theorem}

\proof
 Choose $K>0$ big enough such that 
\begin{equation}\label{bconst}
\mathpzc{b}_m^j:=\bigg(\frac{K \bar{r}_j}{\bar{r}'_{jm}}\frac{s^j_m}{\bar{x}_m}-1\bigg)\frac{1}{\bar{x}_m}>0,
\end{equation}
for any $j$ and $m$. Then consider
\begin{equation}\label{a}
\mathpzc{a}_j := K \bar{r}_j \prod_{m\in\mathbf{M}} \Bigg( \frac{\bar{x}_m}{(1+\mathpzc{b}^j_m \bar{x}_m)}\Bigg)^{-s^j_m}.
\end{equation}
A simple computation shows that the Michaelis-Menten function
$$f_j(x_m):=\mathpzc{a}_j  \prod_{ m}\Bigg( \frac{x_m}{(1+\mathpzc{b}^j_m x_m)}\Bigg)^{s^j_m}$$
satisfies \eqref{mmconstr}.
\endproof
Theorem \ref{MM} shows that any spectral condition of the symbolic Jacobian $G$ can be achieved by a choice of Michaelis-Menten constants. Moreover, it can be satisfied at an equilibrium concentration.
\section{Examples}\label{examplesec}
We present two examples to clarify the results and methods.
\subsection{Example I}
Consider the following reaction network.\\
\begin{minipage}{.6\textwidth}
\begin{align*}       
        A + 2B &\underset{1}{\rightarrow} 2A\\
        \quad\quad\;\;B &\underset{2}{\rightarrow} A+C\\
        \quad\quad\;\;C &\underset{3}{\rightarrow} B
\end{align*}
\end{minipage}%
\begin{minipage}{.1\textwidth}
\begin{align*}
  \quad\quad\;\;A &\underset{0}{\rightarrow} \\
\quad\quad\;\;\quad&\underset{F_B}{\rightarrow}B  
\end{align*}
\end{minipage}

There are only three 3-Child-Selections $\mathbf{J}^{(3)}_1,\mathbf{J}^{(3)}_2,\mathbf{J}^{(3)}_3$.
    \begin{align*}
         \mathbf{J}^{(3)}_1(A,B,C)&=(1,2,3) \quad\quad\quad \alpha_{\mathbf{J}^{(3)}_1}=\operatorname{det}S[\mathbf{J}^{(3)}_1]=\operatorname{det}\begin{pmatrix}
             1 & 1 & 0\\
             -2 & -1 & 1\\
             0 & 1 & -1\\
         \end{pmatrix}=-2;\\
         \mathbf{J}^{(3)}_2(A,B,C)&=
         (0,2,3)\quad\quad\quad \alpha_{\mathbf{J}^{(3)}_2}=\operatorname{det}S[\mathbf{J}^{(3)}_2]=\operatorname{det}\begin{pmatrix}
             -1 & 1 & 0\\
             0 & -1 & 1\\
             0 & 1 & -1\\
         \end{pmatrix}=0;\\
         \mathbf{J}^{(3)}_3(A,B,C)&=
         (0,1,3)\quad\quad\quad \alpha_{\mathbf{J}^{(3)}_3}=\operatorname{det}S[\mathbf{J}^{(3)}_3]=\operatorname{det}\begin{pmatrix}
             -1 & 1 & 0\\
             0 & -2 & 1\\
             0 & 0 & -1\\
         \end{pmatrix}=-2.
    \end{align*}
From $\alpha_{\mathbf{J}^{(3)}_i} \le0$ for $i=1,2,3$, Corollary \ref{corinst2} implies that $a_3=\operatorname{det}G$ is of fixed sign. Since $S[\mathbf{J}^{(3)}_3]$ is a stable matrix, Corollary \ref{stable} implies that the network admits stability. 
From the fact that reaction 1 is autocatalytic, Corollary \ref{autocat} implies that the network admits instability. Thus Corollary \ref{corpi} implies that the Jacobian $G$ admits purely imaginary eigenvalues. In detail, the dynamical system reads:
\begin{equation}\label{msex}
\begin{cases}
    \dot{x}_A=-r_0(x_A) + r_1(x_A,x_B)+r_2(x_B) \\
    \dot{x}_B=-2r_1(x_A,x_B)-r_2(x_B)+r_3(x_C)+F_B\\
    \dot{x}_C=r_2(x_B)-r_3(x_C)
    \end{cases},
\end{equation}
with symbolic Jacobian
$$G(\mathbf{r}')=\begin{pmatrix}
    -r'_{0A}+r'_{1A} & r'_{1B}+r'_{2B} & 0\\
  -2r'_{1A} & -2r'_{1B}-r'_{2B} & r'_{3C}\\
  0 & r'_{2B} & -r'_{3C}
\end{pmatrix}.$$
We note that at $\mathbf{r}'=(r'_{0A},r'_{1A},r'_{1B},r'_{2B},r'_{3C})=(1,1,1,1,1)$ the Jacobian $G$ is still stable. We proceed by making dominant $r'_{1A}$, the `source of instability'. A simple computation shows that for 
\begin{equation}\label{eqpd}
    \bar{\mathbf{r}}'=(\bar{r}'_{0A},\bar{r}'_{1A},\bar{r}'_{1B},\bar{r}'_{2B},\bar{r}'_{3C})=(1,3.5,1,1,1),
    \end{equation}
$G(\bar{\mathbf{r}}')$ has eigenvalues:
\begin{equation}\label{eigv}
    \lambda_1= i\sqrt{6},\quad \lambda_2=-i\sqrt{6},\quad\lambda_3=-1.5.
\end{equation}

We realize this spectrum at an equilibrium of a Michaelis-Menten system, as prescribed by Theorem \ref{MM}. We fix arbitrarily equilibrium concentrations $\bar{x}_A=\bar{x}_B=\bar{x}_C=1$ and sufficiently large equilibrium rates $\bar{r}_1=\bar{r}_2=\bar{r}_3=7$, $\bar{r}_0=\bar{F}_B=14$. Fixing the value of partial derivatives $\bar{\mathbf{r}}'$ as \eqref{eqpd} we obtain the reaction functions
\begin{equation}\label{rf}
f(x)=
\begin{cases}
    f_1(x)=2744 \frac{x_Ax^2_B}{(1+x_A)(1+13x_b)^2}\\
    f_2(x)=49\frac{x_B}{1+6x_B}\\
    f_3(x)=49\frac{x_C}{1+6x_C}\\
    f_0(x)=96\frac{x_A}{1+6x_A}\\
    f_{F_B}\equiv14
\end{cases}.
\end{equation}

The system \eqref{msex} with reaction functions \eqref{rf} at the equilibrium $\bar{x}$ has Jacobian with eigenvalues \eqref{eigv}.


\subsection{Example II}

Here we look for two distinct parameter regions where multistationarity and oscillations appear in numerical simulations. We consider 5 species $A,B,C,D,E$, and 8 reactions.
\begin{minipage}{.6\textwidth}
\begin{align*}       
        A &\underset{1}{\rightarrow} B+C\\
        \quad\quad\;\;B &\underset{2}{\rightarrow} C\\
\quad\quad\;\;C+D&\underset{3}{\rightarrow} A\\
C &\underset{4}{\rightarrow} E     
\end{align*}
\end{minipage}%
\begin{minipage}{.1\textwidth}
\begin{align*}
  \quad\quad\;\;D
  &\underset{5}{\rightarrow} 2B\\
\quad\quad\;\;\quad D+E&\underset{6}{\rightarrow}2E\\
E&\underset{7}{\rightarrow}\\
&\underset{F_D}{\rightarrow} D
\end{align*}
\end{minipage}

The associated system of differential equations is
\begin{equation}\label{example2}
    \begin{cases}
    \dot{x}_A=-r_1(x_A)+r_3(x_C,x_D)\\
    \dot{x}_B=r_1(x_A)-r_2(x_B)+2r_5(x_D)\\
    \dot{x}_C=r_1(x_A)+r_2(x_B)-r_3(x_C,x_D)-r_4(x_C)\\
    \dot{x}_D=-r_3(x_C,x_D)-r_5(x_D)-r_6(x_D,x_E)+F_D\\
    \dot{x}_E=r_4(x_C)+r_6(x_D,x_E)-r_7(x_E)\\
    \end{cases},
\end{equation}
with symbolic Jacobian:
 \begin{equation}\label{symjex2}
 G(\mathbf{r}')=
    \begin{pmatrix}
        -r'_{1A} & 0 & r'_{3C} & r'_{3D} & 0\\
        r'_{1A} & -r'_{2B} & 0 & 2r'_{5D} & 0\\
        r'_{1A} & r'_{2B} & -r'_{3C}-r'_{4C} & -r'_{3D} & 0\\
        0 & 0 & -r'_{3C} &-r'_{3D} -r'_{5D}-r'_{6D} & -r'_{6E}\\
        0 & 0 & r'_{4C} & r'_{6D} & r'_{6E} - r'_{7E}
   \end{pmatrix}.
\end{equation}

\paragraph{Hunting multistationarity:}
Consider the two 5-Child-Selections $\mathbf{J}^{(5)}_1$, $\mathbf{J}^{(5)}_2$:
\begin{align*}
\mathbf{J}^{(5)}_1(A,B,C,D)&=\mathbf{J}^{(5)}_2(A,B,C,D)=(1,2,3,5);\\
   \mathbf{J}^{(5)}_1(E)&=7 \text{ and } \mathbf{J}^{(5)}_2(E) = 6.
\end{align*}
It is easy to see that $$\alpha_{\mathbf{J}^{(5)}_1} \alpha_{\mathbf{J}^{(5)}_2}<0,$$
since
$$\operatorname{det}
\begin{pmatrix}
-1 & 0 & 1 & 0 & 0\\
1 & -1 & 0 & 2 & 0\\
1 & 1 & -1 & 0 & 0\\
0 & 0 & -1 & -1 & 0\\
0 & 0 & 0 & 0 & -1
\end{pmatrix}
=-\operatorname{det}
\begin{pmatrix}
  -1 & 0 & 1 & 0 &0\\
1 & -1 & 0 & 2 & 0\\
1 & 1 & -1 & 0 & 0\\
0 & 0 & -1 & -1 & -1\\
0 & 0 & 0 & 0 & 1
\end{pmatrix}.$$
Hence there are two $5$-Child-Selections of different sign
and Corollary \ref{corinst2} implies that there exists a choice of $\mathbf{r}'$ such that the $a_5=\operatorname{det}G=0$. This may indicate a saddle-node bifurcation and the consequent possibility for multistationarity.

\begin{figure}
    \centering
    \includegraphics[scale=0.8]{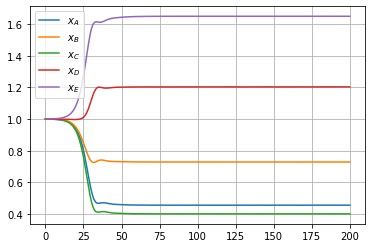}
    \caption{For the system \eqref{example2} with reaction functions \eqref{mmsn}, we have an unstable equilibrium $E_u$ at $x=(1,1,1,1,1)$. In the figure, nearby initial conditions $x(0)=(1.001,1,1,1,1)$ show convergence to the stable equilibrium $E_s\approx(0.45598062, 0.72887142, 0.40193953, 1.20347748, 1.6490945).$}
    \label{fig:multi}
\end{figure}

According to the intuition, we fix the values in the symbolic Jacobian \eqref{symjex2} as follows $$\bar{r}'_{1A}=\bar{r}'_{2B}=\bar{r}'_{3C}=\bar{r}'_{3D}=\bar{r}'_{5D}=\bar{r}'_{6E}=1,\quad\text{and}\quad \bar{r}'_{4C}=\bar{r}'_{6D}=\frac{1}{4}.$$
Then, the Jacobian determinant $$\operatorname{det}G(r'_{7E})=\frac{3}{4} - \frac{21}{16} r'_{7E}$$ changes sign at the bifurcation value $$r'^*_{7E}=\frac{12}{21}.$$
It is straightforward to verify and prove that such a change of stability corresponds to a saddle-node bifurcation. We do not go into this detail, but we simply guess accordingly an area of multistationarity for 
$$\bar{r}'_{7E}:=\frac{1}{2}<\frac{12}{21}.$$ We fix arbitrarily the unstable equilibrium at $\bar{x}=(1,1,1,1,1)$  and equilibrium fluxes $\bar{\mathbf{r}}=(\bar{r}_1,\bar{r}_2,\bar{r}_3,\bar{r}_4,\bar{r}_5,\bar{r}_5,\bar{r}_7,F_D)=(2,4,2,4,1,2,6,5)$. Operating as described in the previous example, we recover the Michaelis-Menten rates:
\begin{equation}\label{mmsn}
\begin{pmatrix}
    r_1(x_A)\\
    r_2(x_B)\\
    r_3(x_C,x_D)\\
    r_4(x_C)\\
    r_5(x_D)\\
    r_6(x_D,x_E)\\
    r_7(x_E)\\
    F_D
\end{pmatrix}=
\begin{pmatrix}
    4\frac{x_A}{1+x_A}\\
    16\frac{x_B}{1+3x_B}\\
    8\frac{x_C}{1+x_C}\frac{x_D}{1+x_D}\\
    64 \frac{x_C}{1+15x_C}\\
     x_D\\
    32 \frac{x_D}{1+7x_D} \frac{x_E}{1+x_E}\\
    72 \frac{x_E}{1+11x_E}\\
    5
\end{pmatrix},
\end{equation}
for which there are two equilibria: one unstable $$E_u=(1,1,1,1,1),$$
set by the construction, and one stable $$E_s\approx(0.45598062, 0.72887142, 0.40193953, 1.20347748, 1.6490945),$$ 
numerically computed. See Figure \ref{fig:multi}.

\paragraph{Hunting oscillations:} Consider the  $4$-Child-Selection $\bar{\mathbf{J}}^{(4)}(A,B,C,D)=(1,2,3,5)$. The 4 eigenvalues of
$$S[\bar{\mathbf{J}}^{(4)}]=
\begin{pmatrix}
-1 & 0 & 1 & 0\\
1 & -1 & 0 & 2\\
1 & 1 & -1 & 0\\
0 & 0 & -1 & -1
\end{pmatrix}
$$
are all stable:
\begin{equation}\label{eigstablehopf}
(\lambda_1,\lambda_2,\lambda_3,\lambda_4)\approx(-2.32472,-1,-0.337641 \pm 0.56228 i).
\end{equation}
Moreover, for the $3$-Child-Selection $\bar{\mathbf{J}}^{(3)}(A,B,C)=(1,2,3)$, the associated matrix
$$S[\bar{\mathbf{J}}^{(3)}]=\begin{pmatrix}
   -1 & 0 & 1\\
1 & -1 & 0\\
1 & 1 & -1 
\end{pmatrix}$$
has determinant 
\begin{equation}\label{detunstablehopf}
    \operatorname{det}S[\bar{\mathbf{J}}^{(3)}]=1=(-1)^{3-1}.
\end{equation} 
Theorem \ref{pithm} concludes the existence of a choice of $\mathbf{r}'$ such that $G(\mathbf{r}')$ has purely imaginary eigenvalues. This indicates the possibility of a Hopf bifurcation and consequently oscillations.

In the same spirit as the hunt for multistationarity, we can find values for which the system shows numerically periodic oscillations. Let us focus on the Cauchy-Binet component associated to $\bar{\mathbf{J}}^{(4)}$:
$$
G[\bar{\mathbf{J}}^{(4)}]=
\begin{pmatrix}
-r'_{1A} & 0 & r'_{3C} & 0\\
r'_{1A} & -r'_{2B} & 0 & 2r'_{5D}\\
r'_{1A} & r'_{2B} & -r'_{3C} & 0\\
0 & 0 & -r'_{3C} & -r'_{5D}
\end{pmatrix}.
$$
Fix $\bar{r}'_{1A}=\bar{r}'_{2B}=\bar{r}'_{3C}=2$. At $r'_{5D}=2$ the matrix is stable because of \eqref{eigstablehopf}, but it loses stability as $r'_{5D} \rightarrow 0$, because of \eqref{detunstablehopf}. Since $G[\bar{\mathbf{J}}^{(4)}]$ is invertible for all $r'_{5D}>0$, the loss of stability must happen in form of purely imaginary eigenvalues crossing, hinting at a Hopf bifurcation.

\begin{figure}
    \centering
\includegraphics[width=\textwidth]{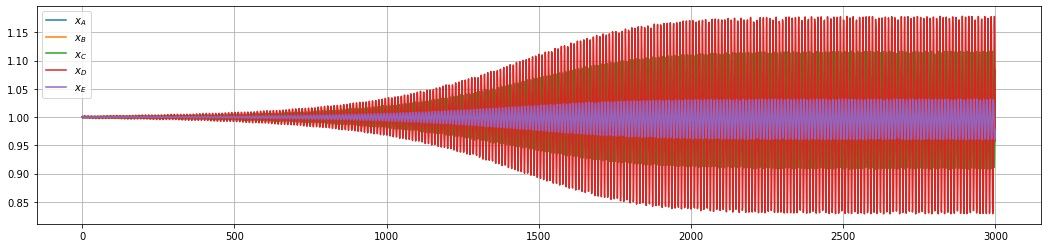}
    \caption{The picture shows numerical simulations for the system \eqref{example2} with reaction functions \eqref{mmhopf}. The system possesses an unstable equilibrium $E_u=(1,1,1,1,1)$. In the plot, nearby initial condition $x(0)=[1.001, 1, 1, 1, 1]$ shows convergence to a stable periodic orbit.} 
    \label{fig:Hopf}
\end{figure}

\begin{figure}
    \centering
    \includegraphics[width=\textwidth]{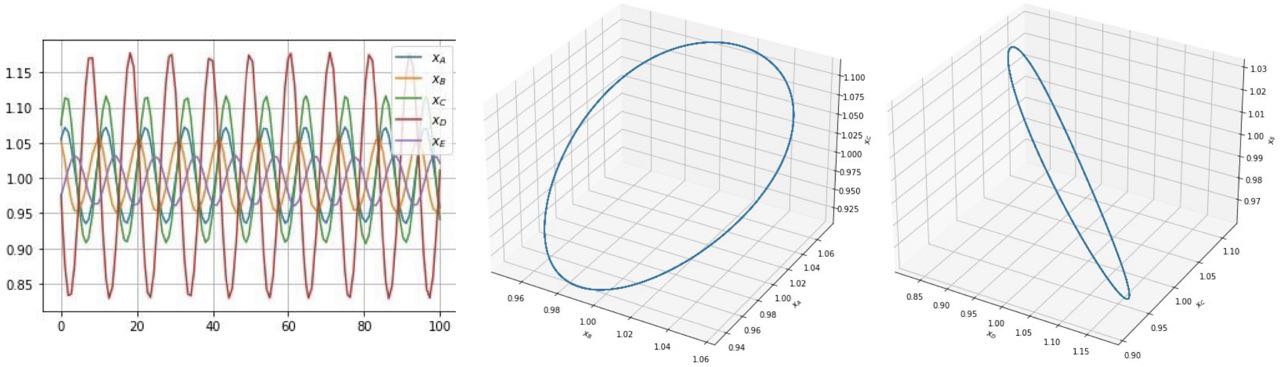}
    \caption{The picture shows numerical simulations for the system \eqref{example2} with reaction functions \eqref{mmhopf}. The initial condition $x(0)=[1.05383533, 1.05121633, 1.07542591, 0.97648847, 0.97416258]$ is chosen in proximity of the periodic orbit. The plot on the left shows the time-evolution of the concentrations $x(t)$, while the plots in the center and on the right depict in 3d the periodic orbit for $(x_A,x_B,x_C)$ and $(x_C,x_D,x_E)$, respectively.}
    \label{fig:Hopf2}
\end{figure}

Imitating the proof of Lemma \ref{piprop}, we fix the other symbols in the Jacobian \eqref{symjex2} small enough as follows 
\begin{align*}
\bar{r}'_{3D}=\bar{r}'_{6E}=\bar{r}'_{7E}&=0.5;\\
\bar{r}'_{4C}&=0.25;\\
\bar{r}'_{6D}&=0.03125.\\
\end{align*}
We use $r'_{5D}$ as a bifurcation parameter. At $r'_{5D}=2$, the Jacobian $G(r'_{5D})$ has 5 negative-real-part eigenvalues
$$ (\lambda_1,\lambda_2,\lambda_3,\lambda_4,\lambda_5)\approx(-5.1256, -1.95591,-0.136865,-0.781435 \pm 0.974403 i).$$ On the other hand for $\bar{r}'_{5D}=0.5$ we find two complex-conjugate eigenvalues with positive real-part:
$$(\lambda_1,\lambda_2,\lambda_3,\lambda_4,\lambda_5)\approx (-0.136034, -3.57546 \pm 0.623139 i, +0.0028519 \pm 0.597922 i ).$$
For such fixed values of $\bar{\mathbf{r}}'$, we fix arbitrarily the unstable equilibrium at $\bar{x}=(1,1,1,1,1)$  and equilibrium fluxes $\bar{\mathbf{r}}=(\bar{r}_1,\bar{r}_2,\bar{r}_3,\bar{r}_4,\bar{r}_5,\bar{r}_5,\bar{r}_7,F_D)=(2,4,2,4,1,2,6,5)$. We recover the Michaelis-Menten rates:
\begin{equation}\label{mmhopf}
\begin{pmatrix}
    r_1(x_A)\\
    r_2(x_B)\\
    r_3(x_C,x_D)\\
    r_4(x_C)\\
    r_5(x_D)\\
    r_6(x_D,x_E)\\
    r_7(x_E)\\
    F_D
\end{pmatrix}=
\begin{pmatrix}
    2x_A\\
    8\frac{x_B}{1+x_B}\\
    8\frac{x_Cx_D}{1+3x_D}\\
    64 \frac{x_C}{1+15x_C}\\
    2 \frac{x_D}{1+x_D}\\
    512 \frac{x_D}{1+63x_D} \frac{x_E}{1+3x_E}\\
    72 \frac{x_E}{1+11x_E}\\
    5
\end{pmatrix}.
\end{equation}
The unstable equilibrium is encircled by a stable periodic orbit. This suggests and confirms indeed the occurrence of a supercritical Hopf bifurcation. See Figures \ref{fig:Hopf} and \ref{fig:Hopf2}. 

\section{Discussion and outlook}\label{discussion}

We presented a symbolic approach to finding bifurcation points in reaction systems with general kinetics. We considered as free positive symbols any partial derivative $r'_{jm}$ of a reaction rate $r_j$ with respect to a concentration of one of its reactant $m$. This way we interpreted and addressed the Jacobian matrix $G$ as a symbolic matrix. The main tool in this framework are the $n$-Child-Selections $\mathbf{J}^{(n)}$: injective maps that associate species $m$ to reactions $j=\mathbf{J}^{(n)}(m)$ where $m$ participates as reactant and, in particular, $r'_{\mathbf{J}^{(n)}(m)m}\neq0$. We employed a Cauchy-Binet analysis to structurally express any coefficient of the characteristic polynomial of $G$ in terms of Child-Selections (Theorem \ref{CBTHM}). In particular, we gave necessary and sufficient conditions such that each of the coefficients admits a change of sign for some choice of symbols (Corollary \ref{corinst2}). The conditions are in terms of signs of stoichiometric submatrices identified by the Child-Selections. The $M^{th}$ constant coefficient $a_M$ is the determinant of the Jacobian, whose sign changes hint at bifurcation behavior. More generally, even in presence of conserved quantities, a sign change of the highest non-identically-zero coefficient $a_n$ is a necessary spectral condition to obtain codimension 1 bifurcations such as saddle-node bifurcations, transcritical bifurcations, pitchfork bifurcations, or even codimension 2 bifurcations such as Takens-Bogdanov. Moreover, we analyzed how the stability properties of any Child-Selection $\mathbf{J}^{(n)}$, expressed by the inertia of its associated Cauchy-Binet component $G[\mathbf{J}^{(n)}]$, are inherited by the full Jacobian via a perturbation argument (Theorem \ref{stability}). In our symbolic setting, inheritance is a quite natural property, in contrast with more restrictive settings addressed in the literature 
\cite{Ba18,BaPa18}. Finally, we focused on sufficient conditions that guarantee that the Jacobian $G$ admits purely imaginary eigenvalues, hinting at Hopf bifurcation and consequent oscillatory behavior (Lemma \ref{piprop} and Theorem \ref{pithm}). We underline however that our results do not automatically imply nonresonance conditions such as the absence of other eigenvalues with zero-real part.

There is full validity of the presented symbolic approach within the defined class of monotone chemical functions (Definition \ref{monchemfunct}). Fix any choice of symbols $\bar{\mathbf{r}}'$: `full validity' simply means that we can always find a positive equilibrium $\bar{x}>0$, and monotone chemical reaction functions $f$ such that the Jacobian $G$ of $\dot{x}=Sf(x)$, evaluated at  $\bar{x}$, reads as $G(\bar{\mathbf{r}}')$. More concretely, we noted that the standard Michaelis-Menten kinetics already allows enough parametric freedom to obtain such spectral realization at any chosen equilibrium $\bar{x}$ (Theorem \ref{MM}). The extension to richer classes of functions that contain Michaelis-Menten as a subclass is of course straightforward. This is for example the case of Hill kinetics. On the contrary, for the more restrictive class of mass action kinetics, our symbolic analysis does not directly conclude the realizability of a spectral property at an equilibrium of a mass action system. However, exclusion results that assert the absence of a certain bifurcation behavior are instead fully valid, since any mass action realization of the Jacobian fall within the symbolic realization. Different types of kinetics in the class of monotone chemical functions can be further checked, along similar lines. 

Besides linear conditions, bifurcations theorems require further nondegeneracy conditions involving higher derivatives. This is in no way trivial and it has still to be discussed case-by-case. For existence results of bifurcations along the symbolic lines of this paper see \cite{V23} for an account on saddle-node bifurcations. In a related setting, see Fiedler \cite{F19} who proved \emph{global} Hopf bifurcation \cite{AlexYorke78}  based on some cycle conditions on the symbolic Jacobian. 

For readers who are more interested in finding bifurcations in numerical simulations rather than proving them analytically, our contribution already suggests various possibilities to guess the parameter areas of interest. We exemplified such possibilities in the two toy models of Section \ref{examplesec}: Example I and Example II. {We underline that we do not rely on any Hurwitz computation to conclude the occurrence of the bifurcations. The development of systematic algorithms to detect the spectral conditions for bifurcations falls beyond the scope of the present work, and will be addressed in future research. Suffice it to say that for networks of sufficiently small size we can resort to computer algebra tools  for the computation of the characteristic polynomial. As a consequence, e.g., we can automatically identify $n$-Child-Selections that leads to instability, according to Corollary \ref{cor1}: it is enough to check monomials that carry the `wrong' sign. However, such an approach is doomed to fail for reasonably large and realistic metabolic networks where dozens, if not hundreds, of species and reactions are involved. In this case, non-trivial heuristic algorithms need to be developed for the purpose. }

In conclusion, we point to two open questions:
\begin{enumerate}
\item Which sufficient conditions on an invertible $S[\mathbf{J}^{(n)}]$ guarantee that 
\begin{equation}\label{CBnohopf}
    \operatorname{inertia}G[\mathbf{J}^{(n)}](\mathbf{r}')\equiv \operatorname{inertia}S[\mathbf{J}^{(n)}],
\end{equation}
for any choice of symbols $\mathbf{r}'$?
\item When and how the set of the inertias of all Cauchy-Binet components is sufficient to determine the inertia of $G$?
\end{enumerate}
The first question has been already formulated in the literature in general linear algebra setting \cite{Ku19}. Consider any $M\times M$ real invertible matrix $A$ and a positive symbolic diagonal $M\times M$ matrix $D$. $A$ is said to be \emph{$D$-hyperbolic} if the inertia of $A$ corresponds to the inertia of $AD$ for any choice of $D$. Since $A$ is invertible, $D$-hyperbolicity is equivalent to the absence of purely imaginary eigenvalues of $AD$ for any $D$. For instance, triangular matrices transfer their inertia to $AD$, for any choice of positive $D$, but less trivial classes have to be investigated. The special case in which $A$ is a stable matrix has been addressed more extensively in the literature, see \cite{Gio15} for a focused review. Even in such first simplest case however, a full characterization of $D$-stability remained elusive. The second question is intimately connected to the current network setting. For example, assume that all Cauchy-Binet components $G[\mathbf{J}^{(n)}]$ are stable
 $$\operatorname{inertia}G[\mathbf{J}^{(n)}](\mathbf{r}')=(n,0,0),$$
for all choice of symbols $\mathbf{r}'$. In particular, \eqref{CBnohopf} holds for all Cauchy-Binet components $G[\mathbf{J}^{(n)}](\mathbf{r}')$. We ask which network conditions are sufficient to conclude that $G$ cannot have purely imaginary eigenvalues for the Jacobian $G$ for all choices of symbols, and it is always a stable matrix. A better understanding of these questions would help in clarifying or excluding bifurcation behaviors in chemical reaction networks.

\section{Proofs of Section \ref{main}}\label{proofs}

\subsection*{Proofs of subsection \ref{CB}}

\proof[Proof of Lemma \ref{corinst}]
Via Cauchy-Binet formula, we compute.
\begin{align*}
 \operatorname{det}G^{\mathbf{M}^{(n)}}_{\mathbf{M}^{(n)}} = \operatorname{det}  (SR)^{\mathbf{M}^{(n)}}_{\mathbf{M}^{(n)}}
     = \operatorname{det}  S_{\mathbf{M}^{(n)}} R^{\mathbf{M}^{(n)}}
     = \sum_{|\vartheta|=n}  \operatorname{det} S^\vartheta_{\mathbf{M}^{(n)}} \operatorname{det} R^{\mathbf{M}^{(n)}}_\vartheta.
\end{align*}
We note that $\operatorname{det}R^{\mathbf{M}^{(n)}}_\vartheta \neq 0$ if and only if $\vartheta=\mathbf{J}^{(n)}(\mathbf{M}^{(n)})$, for some $n$-Child-Selection $\mathbf{J}^{(n)}$. In particular, for $\vartheta=\mathbf{J}^{(n)}(\mathbf{M}^{(n)})$ we have
$$\operatorname{det}R^{\mathbf{M}^{(n)}}_\vartheta=\operatorname{sgn}(\mathbf{J}^{(n)}) \prod_{m \in \mathbf{M}^{(n)}} r'_{\mathbf{J}^{(n)}(m)m},$$
where $\operatorname{sgn}(\mathbf{J}^{(n)})$ indicates the \emph{signature} of $\mathbf{J}^{(n)}$ as a map from an ordered set of $n$ elements to another ordered set of $n$ elements. Moreover, for such a choice of $\vartheta$,
$$ \operatorname{det}S^{\vartheta}_{\mathbf{M}^{(n)}}  \operatorname{sgn}(\mathbf{J}^{(n)})=\operatorname{det}S[\mathbf{J}^{(n)}].$$
This yields:
\begin{align}\label{forcor}
\begin{split}
\operatorname{det}G^{\mathbf{M}^{(n)}}_{\mathbf{M}^{(n)}}&=\sum_{\vartheta=\mathbf{J}^{(n)}(\mathbf{M}^{(n)})}  \operatorname{det} S^{\vartheta}_{\mathbf{M}^{(n)}}  \operatorname{sgn}(\mathbf{J}^{(n)}) \prod_{m \in \mathbf{M}^{(n)}} r_{\mathbf{J}^{(n)}(m)m}\\&=\sum_{\mathbf{J}^{(n)}:\mathbf{M}^{(n)}\mapsto \mathbf{E}}  \operatorname{det} S[\mathbf{J}^{(n)}]\operatorname{det}R[\mathbf{J}^{(n)}]\\&=\sum_{\mathbf{J}^{(n)}:\mathbf{M}^{(n)}\mapsto \mathbf{E}} \operatorname{det}G[\mathbf{J}^{(n)}].
\end{split}
\end{align}
\endproof

\proof[Proof of Theorem \ref{CBTHM}]
Any coefficient $a_n$ in the characteristic polynomial $\mathpzc{g}(\lambda)$ of $G$ is the sum of all the principal minors of $G$ of order $n$. We compute.
\begin{align*}
    a_n &= \sum_{\mathbf{M}^{(n)}} \operatorname{det} G^{\mathbf{M}^{(n)}}_{\mathbf{M}^{(n)}}\\
    &= \sum_{\mathbf{M}^{(n)}} \;\; \sum_{\mathbf{J}^{(n)}:\mathbf{M}^{(n)}\mapsto \mathbf{E}}  \operatorname{det}G[\mathbf{J}^{(n)}]\\
    &= \sum_{\mathbf{J}^{(n)}} \operatorname{det}G[\mathbf{J}^{(n)}],
\end{align*}
where the last sum runs on all $n$-Child-Selections.
\endproof

\proof[Proof of Corollary \ref{corinst2}]
 Lemma \ref{corinst} and Theorem \ref{CBTHM} imply that any principal minor $G^{\mathbf{M}^{(n)}}_{\mathbf{M}^{(n)}}$ and any coefficient $a_n$ is a multilinear $n$-homogenous polynomial in the symbols $\mathbf{r}'$. More specifically, the linear monomial summands are of the form 
 $$\alpha_{\mathbf{J}^{(n)}} \prod_{m\in\mathbf{M}^{(n)}} r'_{\mathbf{J}^{(n)}(m)m}.$$
For any multilinear $n$-homogenous polynomials $p[y]$ in the positive variables $y\in \mathbb{R}_{>0}^N$, $$p[y]=\sum_i p_i[y]=\sum_i \alpha_i \prod_{i=1}^n y_i,$$ 
it holds that $p[y]$ is of fixed sign if and only if there are no two monomial summands $p_h$ and $p_l$ that have coefficients with 
$$\alpha_h \; \alpha_l <0.$$
In fact, one direction is trivial: if all nonzero monomial summands have the same nonnegative (resp. nonpositive) sign then the polynomial $p[y]$ is either negative (resp. positive) or zero if $p[y]\equiv0$. Indirectly assume now that there are two nonzero summands of different sign: $p_h[y]>0$ and $p_l[y]<0$. For the value $y_{h(1)}=y_{h(2)}=...=y_{h(n)}=H$ and big enough $H$, $p[y]$ attains positive values, while for the $y_{l(1)}=y_{l(2)}=...=y_{l(n)}=L$ and big enough $L$, $p[y]$ attains negative values. This contradicts the assumption of fixed sign.
\endproof

\proof[Proof of Corollary \ref{pmatrix}]
The Jacobian $G(\mathbf{r}')$ is a $P^{(-)}_0$ matrix for any choice of $\mathbf{r}'$ if and only if all principal minors are of fixed sign $(-1)^{n}$. All principal minors are of fixed sign $(-1)^{n}$ if and only if for all $n$-Child-Selection $\mathbf{J}^{(n)}$ it holds
$$\operatorname{sign}\alpha_{\mathbf{J}^{(n)}}=(-1)^n.$$
\endproof

\subsection*{Proofs of subsection \ref{instabilities}}

\proof[Proof of Theorem \ref{stability}]
Consider the choice of symbols $\mathbf{r}'$ and the following rescaling:
\begin{equation}
    \rho'_{jm}(\varepsilon)=
    \begin{cases}
    \varepsilon r'_{jm} \quad\;\; \text{if $(j,m)\neq (\mathbf{J}^{(n)}(m),m$});\\
    r'_{jm } \quad\quad \text{otherwise}.
    \end{cases}
\end{equation}
 Again, the bold $\boldsymbol{\rho}'(\varepsilon)$ indicates the whole sets of the rescaled symbols $\{\rho'_{jm}(\varepsilon)\}$. The eigenvalues of $G$ do not depend on the labeling of the network. Thus, without loss of generalities, we can consider $\mathbf{M}^{(n)}=\{m_1,...,m_n\}$. For $\varepsilon=0$, the Jacobian $G$ takes the block form:
$$G(\boldsymbol{\rho}'(0))=
\begin{pmatrix}
G[\mathbf{J}^{(n)}] &0\\
... &0
\end{pmatrix}.$$
Now choose $\mathbf{r}'=\mathbf{r}'_1$, then $G(\boldsymbol{\rho}'_1(0))$ has inertia $(\sigma^-_{\mathbf{J}^{(n)}},\sigma^+_{\mathbf{J}^{(n)}},\sigma^0_{\mathbf{J}^{(n)}}+M-n)$. By continuity of the eigenvalues with respect to the matrix entries, there exists $\varepsilon>0$ such that $\operatorname{inertia}(G(\boldsymbol{\rho}'(\varepsilon))=(\sigma^-_G,\sigma^+_G,\sigma^0_G),$
    with $\sigma^+_G \ge \sigma^+_{\mathbf{J}^{(n)}}$, and $\sigma^-_G \ge \sigma^-_{\mathbf{J}^{(n)}}$.
\endproof

\proof[Proof of Corollary \ref{stable}]
By Theorem \ref{stability}, there exists a choice of $\mathbf{r}'$ such that the Jacobian $G(\mathbf{r}')$ has at least $M$ stable eigenvalues. Since $M$ is the entire number of the eigenvalues of $G$, there exists a choice $G(\mathbf{r}')$ for which $G$ is a stable matrix. If $S[\mathbf{J}^{(M)}]$ is a stable matrix, the choice of symbols $r'_{jm} \equiv 1$ makes $G[\mathbf{J}^{(M)}]=S[\mathbf{J}^{(M)}]$ a stable matrix and we can repeat the argument.
\endproof

\proof[Proof of Corollary \ref{cor1}]
If $\operatorname{sign}\alpha_{\mathbf{J}^{(n)}}=(-1)^{n-1}$, then the parity of the eigenvalues with positive real part is odd and in particular there is at least one: $\sigma^+_{\mathbf{J}^{(n)}}>0$. Theorem \ref{stability} implies that there exists a choice of $\mathbf{r}'$ such that $\sigma^+_G \ge \sigma^+_{\mathbf{J}^{(n)}}>0$, and thus the network admits instability.
\endproof

\proof[Proof of Corollary \ref{exclusive}]
Corollaries \ref{pmatrix} and \ref{cor1} directly imply this.
\endproof

\proof[Proof of Corollary \ref{autocat}]
Consider the 1-Child-Selection $\bar{\mathbf{J}}^{(1)}$ on the set $\bar{\mathbf{M}}^{(1)}=\{ m \}$ such that 
$\bar{\mathbf{J}}^{(1)}(m)=j_{aut}$. The coefficient behavior of $\bar{\mathbf{J}}^{(1)}$ is
$$\bar{s}^{j_{aut}}_m-s^{j_{aut}}_m>0,$$
with $\operatorname{sign}\alpha_{\bar{\mathbf{J}}^{(1)}}=(-1)^{1-1}$. Corollary \ref{cor1} implies that the network admits instability.
\endproof

\subsection*{Proofs of subsection \ref{purelyim}}

\proof[Proof of Lemma \ref{piprop}]
The existence of a choice of symbols $\mathbf{r}'_1$ such that $G[\mathbf{J}^{(n)}]$ is stable implies that the associated integer matrix $S[\mathbf{J}^{(n)}]$ is invertible. In fact, for $\mathbf{r}'_1$,
$$0\neq(-1)^n=\operatorname{sign}\operatorname{det}G[\mathbf{J}^{(n)}](\mathbf{r}'_1)=\operatorname{sign}\operatorname{det}S[\mathbf{J}^{(n)}]\operatorname{sign}\operatorname{det}R[\mathbf{J}^{(n)}](\mathbf{r}_1')=\operatorname{sign}\operatorname{det}S[\mathbf{J}^{(n)}].$$
Moreover, the same equalities hold at any $\mathbf{r}'$, which implies that $G[\mathbf{J}^{(n)}]$ is invertible for any choice of symbols.

Consider now any continuous curve in symbol space
$$\gamma(\mu):\quad [0,1] \quad \mapsto \quad \mathbb{R}^{n},$$
joining $\mathbf{r}'_1$ and $\mathbf{r}'_2$. That is, 
$$\begin{cases}
    \gamma(0)=\mathbf{r}'_1\\
    \gamma(1)=\mathbf{r}'_2.
\end{cases}$$
The intermediate value theorem implies that there exists at least one value $\mu^*$ where the real part of at least one eigenvalue of $G[\mathbf{J}^{(n)}]$ changes sign. In particular,
$$\operatorname{inertia}G[\mathbf{J}^{(n)}](\gamma(\mu^*))=(\sigma^-,\sigma^+,\sigma^0),$$
with $\sigma^0\ge1$. Invertibility of $G[\mathbf{J}^{(n)}]$ for any $\mathbf{r}'$ excludes crossings by real eigenvalues zero, thus implying a purely imaginary crossing by (at least one) pair of eigenvalues. Therefore,  $G[\mathbf{J}^{(n)}](\gamma(\mu^*))$ possesses purely imaginary eigenvalues. Note also that $G[\mathbf{J}^{(n)}](\gamma(\mu))$ identifies a curve $\Lambda(\mu)$
$$\Lambda:\quad [0,1] \quad \mapsto \quad \mathbb{C}^{n},$$
such that each coordinate $\Lambda_i$ correspond to the curve of the eigenvalue $\lambda_i$. Without loss of generalities, we consider that the purely imaginary crossing at $\mu^*$ happens for the eigenvalues $\lambda_1$ and $\lambda_2$.

Now we consider again the rescale of the symbols $\mathbf{r}'$:
\begin{equation}\label{rescale2}
    \rho'_{jm}(\varepsilon)=
    \begin{cases}
    \varepsilon r'_{jm} \quad\;\; \text{if $(j,m)\neq (\mathbf{J}^{(n)}(m),m$});\\
    r'_{jm } \quad\quad \text{otherwise},
    \end{cases}
\end{equation}
and let the bold $\boldsymbol{\rho}'(\varepsilon)$ indicate the whole sets of the rescaled symbols $\{\rho'_{jm}(\varepsilon)\}$. Without loss of generalities in the labeling of the network, the Jacobian $G$ reads at $\varepsilon=0$ as a block matrix
$$G(\boldsymbol{\rho}'(0))=
\begin{pmatrix}
G[\mathbf{J}^{(n)}] &0\\
... &0
\end{pmatrix}.$$
Fix any value $\bar{r}'_{jm}$ for $(j,m)\neq (\mathbf{J}^{(n)}(m),m)$. We can still consider the curves $\gamma(\mu)$ and $\Lambda(\mu)$, just as before. Of course, now $\Lambda(\mu)$ does not identify the curve of the entire spectrum of $G$, but only the curve of $n$ eigenvalues out of $M$. Let $\gamma^\varepsilon(\mu)$ be the continuous curve obtained by perturbing $\gamma$ according to the rescaling \eqref{rescale2}, and $\Lambda^\varepsilon(\mu)$ the associated curve in the eigenvalue space. For $\varepsilon$ small enough, the following three observations hold:
\begin{enumerate}
\item $\Lambda^\varepsilon(0)$ identifies $n$ stable eigenvalues;
\item $\Lambda^\varepsilon(1)$ identifies at least $2$ unstable eigenvalues;
\item $\Lambda^\varepsilon(\mu) \neq 0$ for all $\mu$. 
\end{enumerate}
Again by the intermediate value theorem, there exists a value $\mu^{**}$ such that $\Lambda^\varepsilon(\mu^{**})$ identifies at least one pair of purely imaginary eigenvalues. Hence, $G(\gamma^\varepsilon(\mu^{**}))$ possesses purely imaginary eigenvalues. 
\endproof

\proof[Proof of Corollary \ref{corlem}]
Since $S[\mathbf{J}^{(n)}]$ is stable, for $r'_{jm}\equiv1$ we get that
$G[\mathbf{J}^{(n)}]=S[\mathbf{J}^{(n)}]$
is stable. On the other hand, consider the network $(\mathbf{M}^{(n)},\mathbf{J}^{(n)}(\mathbf{M}^{(n)}))$ with stoichiometry restricted only to the species $m\in \mathbf{M}^{(n)}$. Corollary \ref{exclusive} implies that there exists a choice of $\mathbf{r}'$ such that $G[\mathbf{J}^{(n)}]$ is unstable. Lemma \ref{piprop} concludes the statement.
\endproof

\proof[Proof of Theorem \ref{pithm}]
We consider again a rescale of $\mathbf{r}'$:
\begin{equation}
    \rho_{jm}(\varepsilon):=
    \begin{cases}
    \varepsilon r_{jm} \quad\;\; \text{if $(j,m)\neq (\mathbf{J}^{(\bar{n})}(m),m$) with $\mathbf{J}^{(\bar{n})}\in\mathpzc{C}$};\\
    r_{jm } \quad\quad \text{otherwise},
    \end{cases}
\end{equation}
and let the bold $\boldsymbol{\rho}'(\varepsilon)$ indicate the whole sets of the rescaled symbols $\{\rho'_{jm}(\varepsilon)\}$. We define the $\bar{n}\times\bar{n}$ matrix $G[\mathpzc{C}]$ as 
$$G[\mathpzc{C}]:=G(\boldsymbol{\rho'}(0))^{\bar{\mathbf{M}}^{(n)}}_{\bar{\mathbf{M}}^{(n)}},$$
that is, $G[\mathpzc{C}]$ is the $\bar{n}$-principal minor of $G(\boldsymbol{\rho}'(0))$ considering only rows and columns associated to the species in $\bar{\mathbf{M}}^{(n)}$. Since $S[\mathbf{J}_1^{(\bar{n})}]$ is a stable matrix, $\operatorname{sign}\alpha_{\mathbf{J}_1^{(\bar{n})}}=(-1)^{\bar{n}}$, and hence $G[\mathpzc{C}]$ is invertible for any choice $\boldsymbol{\rho}'(0)$. Indeed, by construction,
$$\operatorname{sign}\operatorname{det}G[\mathpzc{C}]=\operatorname{sign}\bigg(\sum_{\mathbf{J}^{(\bar{n})}\in \mathpzc{C}} \operatorname{det} G[\mathbf{J}^{(\bar{n})}]\bigg)=(-1)^{\bar{n}}.$$
Moreover, assumption 1 also implies that there exists $\boldsymbol{\rho}'_1(0)$ such that $G[\mathpzc{C}](\boldsymbol{\rho}'_1(0))$ is stable, and assumption 2 implies that there exists $\boldsymbol{\rho}'_2(0)$ such that $G[\mathpzc{C}](\boldsymbol{\rho}'_2(0))$ is unstable, in same spirit as in Theorem \ref{stability}.

The statement follows now in total analogy as Lemma \ref{piprop}, and we just sketch it for self-consistency of the present proof. We have that $G[\mathpzc{C}]$ is invertible for any choice of symbols and there exist two choices $\boldsymbol{\rho}_1'(0),\boldsymbol{\rho}_2'(0)$ of symbols such that $G[\mathpzc{C}]$ is stable and unstable, respectively. Therefore we can find a choice of symbols for which $G[\mathpzc{C}]$ changes stability at purely imaginary eigenvalues, via intermediate value theorem. For $\varepsilon$ small enough the same argument holds for $G$, as the curve connecting $\boldsymbol{\rho}_1'(0)$ and $\boldsymbol{\rho}_2'(0)$ is perturbed in an open region away from zero.
\endproof

\proof[Proof of Corollary \ref{corpi}]
As in proofs of Lemma \ref{piprop} and Theorem \ref{pithm}, consider any continuous curve $\gamma(\mu)$ in symbol space such that at $\gamma(0)$ the Jacobian $G$ is stable and at $\gamma(1)$ the Jacobian $G$ is unstable. The intermediate value theorem implies a loss of stability along the curve $\gamma$. Since the Jacobian is always invertible, the loss of stability happens as purely imaginary eigenvalues crossing.
\endproof

\bibliography{references.bib}

\providecommand{\bysame}{\leavevmode\hbox to3em{\hrulefill}\thinspace}
\providecommand{\MR}{\relax\ifhmode\unskip\space\fi MR }
\providecommand{\MRhref}[2]{%
  \href{http://www.ams.org/mathscinet-getitem?mr=#1}{#2}
}
\providecommand{\href}[2]{#2}
\begin{thebibliography}{10}

\bibitem{AlexYorke78}
James~C. Alexander and James~A. Yorke, \emph{Global bifurcations of periodic
  orbits}, American Journal of Mathematics \textbf{100} (1978), no.~2,
  263--292.

\bibitem{AnBaPa13}
David Angeli, Murad Banaji, and Casian Pantea, \emph{Combinatorial approaches
  to hopf bifurcations in systems of interacting elements}, Communications in
  Mathematical Sciences \textbf{12} (2013), no.~6, 1101--1133.

\bibitem{Ang07}
David Angeli, Patrick De~Leenheer, and Eduardo~D Sontag, \emph{A {P}etri net
  approach to the study of persistence in chemical reaction networks},
  Mathematical biosciences \textbf{210} (2007), no.~2, 598--618.

\bibitem{Ba18}
Murad Banaji, \emph{Inheritance of oscillation in chemical reaction networks},
  Applied Mathematics and Computation \textbf{325} (2018), 191--209.

\bibitem{BaBo23}
Murad Banaji and Bal{\'a}zs Boros, \emph{The smallest bimolecular mass action
  reaction networks admitting andronov--hopf bifurcation}, Nonlinearity
  \textbf{36} (2023), no.~2, 1398.

\bibitem{Ba-07}
Murad Banaji, Pete Donnell, and Stephen Baigent, \emph{P matrix properties,
  injectivity, and stability in chemical reaction systems}, SIAM Journal on
  Applied Mathematics \textbf{67} (2007), no.~6, 1523--1547.

\bibitem{BaPa18}
Murad Banaji and Casian Pantea, \emph{The inheritance of nondegenerate
  multistationarity in chemical reaction networks}, SIAM Journal on Applied
  Mathematics \textbf{78} (2018), no.~2, 1105--1130.

\bibitem{Boros21}
Bal{\'a}zs Boros and Josef Hofbauer, \emph{Oscillations in planar
  deficiency-one mass-action systems}, Journal of Dynamics and Differential
  Equations (2021), 1--23.

\bibitem{BruSha09}
Richard~A. Brualdi and Bryan~L. Shader, \emph{Matrices of sign-solvable linear
  systems}, Cambridge University Press, 2009.

\bibitem{CarstenHopfExclusion19}
Carsten Conradi, Elisenda Feliu, and Maya Mincheva, \emph{On the existence of
  hopf bifurcations in the sequential and distributive double phosphorylation
  cycle}, Mathematical biosciences and engineering: MBE \textbf{17} (2019),
  no.~1, 494--513.

\bibitem{Conradi2007}
Carsten Conradi, Dietrich Flockerzi, and Jorg Raisch, \emph{Saddle-node
  bifurcations in biochemical reaction networks with mass action kinetics and
  application to a double-phosphorylation mechanism}, 2007 American control
  conference, IEEE, 2007, pp.~6103--6109.

\bibitem{conradietal19}
Carsten Conradi, Maya Mincheva, and Anne Shiu, \emph{Emergence of oscillations
  in a mixed-mechanism phosphorylation system}, Bulletin of mathematical
  biology \textbf{81} (2019), no.~6, 1829--1852.

\bibitem{DomKirk09}
Mirela Domijan and Markus Kirkilionis, \emph{Bistability and oscillations in
  chemical reaction networks}, Journal of Mathematical Biology \textbf{59}
  (2009), no.~4, 467--501.

\bibitem{Errami2015}
Hassan Errami, Markus Eiswirth, Dima Grigoriev, Werner~M Seiler, Thomas Sturm,
  and Andreas Weber, \emph{Detection of hopf bifurcations in chemical reaction
  networks using convex coordinates}, Journal of Computational Physics
  \textbf{291} (2015), 279--302.

\bibitem{F19}
Bernold Fiedler, \emph{Global {H}opf bifurcation in networks with fast feedback
  cycles}, Discrete and Continuous Dynamical Systems - S \textbf{0} (2020),
  no.~1937-1632\_2019\_0\_144.

\bibitem{Gat2005}
Karin Gatermann, Markus Eiswirth, and Anke Sensse, \emph{Toric ideals and graph
  theory to analyze {H}opf bifurcations in mass action systems}, Journal of
  Symbolic Computation \textbf{40} (2005), no.~6, 1361--1382.

\bibitem{Gio15}
Giorgio Giorgi and Cesare Zuccotti, \emph{An overview on d-stable matrices},
  Department of Economics and Management DEM Working Paper Series (2015).

\bibitem{GuHo84}
John Guckenheimer and Philip Holmes, \emph{Nonlinear oscillations, dynamical
  systems and bifurcations of vector fields}, Springer, 1984.

\bibitem{Hell16}
Juliette Hell and Alan~D. Rendall, \emph{Sustained oscillations in the map
  kinase cascade}, Mathematical Biosciences \textbf{282} (2016), 162--173.

\bibitem{Hersh92}
Daniel Hershkowitz, \emph{Recent directions in matrix stability}, Linear
  Algebra and its Applications \textbf{171} (1992), 161--186.

\bibitem{Hess71}
Benno Hess and Arnold Boiteux, \emph{Oscillatory phenomena in biochemistry},
  Annual review of biochemistry \textbf{40} (1971), no.~1, 237--258.

\bibitem{Hill10}
Archibald~Vivian Hill, \emph{The possible effects of the aggregation of the
  molecules of haemoglobin on its dissociation curves}, J. Physiol. \textbf{40}
  (1910), 4--7.

\bibitem{Holling65}
Crawford~Stanley Holling, \emph{The functional response of predators to prey
  density and its role in mimicry and population regulation}, The Memoirs of
  the Entomological Society of Canada \textbf{97} (1965), no.~S45, 5--60.

\bibitem{HFJ72}
Fritz Horn and Roy Jackson, \emph{General mass action kinetics}, Archive for
  Rational Mechanics and Analysis \textbf{47} (1972), no.~2, 81--116.

\bibitem{JKD77}
Clark Jeffries, Victor Klee, and Pauline Van~den Driessche, \emph{When is a
  matrix sign stable?}, Canadian Journal of Mathematics \textbf{29} (1977),
  no.~2, 315--326.

\bibitem{Knill14}
Oliver Knill, \emph{Cauchy--binet for pseudo-determinants}, Linear Algebra and
  its Applications \textbf{459} (2014), 522--547.

\bibitem{Ku19}
Olga~Y. Kushel, \emph{Unifying matrix stability concepts with a view to
  applications}, SIAM Review \textbf{61} (2019), no.~4, 643--729.

\bibitem{LOD18}
Jephian C.-H. Lin, Dale~D. Olesky, and Pauline van~den Driessche, \emph{Sign
  patterns requiring a unique inertia}, Linear Algebra and its Applications
  \textbf{546} (2018), 67--85.

\bibitem{Liu94}
Wei-Min Liu, \emph{Criterion of hopf bifurcations without using eigenvalues},
  Journal of Mathematical Analysis and Applications \textbf{182} (1994), no.~1,
  250--256.

\bibitem{MM13}
L.~Michaelis and M.~L. Menten, \emph{Die kinetik der invertinwirkung}, Biochem.
  Z. \textbf{49} (1913), 333--369.

\bibitem{Monod49}
Jacques Monod, \emph{The growth of bacterial cultures}, Annual review of
  microbiology \textbf{3} (1949), no.~1, 371--394.

\bibitem{ThomKauf01}
Ren{\'e} Thomas and Marcelle Kaufman, \emph{Multistationarity, the basis of
  cell differentiation and memory. {I}. {S}tructural conditions of
  multistationarity and other nontrivial behavior}, Chaos: An Interdisciplinary
  Journal of Nonlinear Science \textbf{11} (2001), no.~1, 170--179.

\bibitem{VGB20}
Nicola Vassena, \emph{Good and bad children in metabolic networks},
  Mathematical Biosciences and Engineering \textbf{17} (2020), no.~6,
  7621--7644.

\bibitem{V22}
\bysame, \emph{Structural obstruction to the simplicity of the eigenvalue zero
  in chemical reaction networks}, arXiv preprint arXiv:2205.12655 (2022).

\bibitem{V23}
\bysame, \emph{Structural conditions for saddle-node bifurcations in chemical
  reaction networks}, SIAM Journal on Applied Dynamical Systems \textbf{22}
  (2023), no.~3, 1639--1672.

\bibitem{VasStad23}
Nicola Vassena and Peter~F Stadler, \emph{Unstable cores are the source of
  instability in chemical reaction networks}, arXiv preprint arXiv:2308.11486
  (2023).

\end{thebibliography}
\bibliographystyle{amsplain}

\end{document}